\documentclass[12pt] {article}
\usepackage{amsmath,amsthm}
\usepackage{amssymb,latexsym}
\usepackage{epic,eepic}
\usepackage{amscd}
\usepackage{epic,eepic}
\title{Pluripotential theory on quaternionic manifolds.}
\date{}
\author{Semyon Alesker \footnote{Partially supported by ISF grant 701/08.}
\\  { \normalsize Department of Mathematics, Tel Aviv University, Ramat Aviv}
\\  { \normalsize 69978 Tel Aviv, Israel }
\\ {\normalsize e-mail: semyon@post.tau.ac.il}}

\def\alp{\alpha}
\def\ome{\omega}
\def\Ome{\Omega}

\def\to{\rightarrow}
\def\qed { Q.E.D. }

\def\cct{{}\!^ {\mathbb{C}} T}
\def\cchh{{}\!^ {\mathbb{C}} \mathbb{H}}

\def\ccb{{}\!^ {\mathbb{C}} \mathcal{B}}

\def\RR{\mathbb{R}}
\def\CC{\mathbb{C}}

\def\ZZ{\mathbb{Z}}
\def\HH{\mathbb{H}}

\def\PP{\mathbb{P}}

\swapnumbers
\newtheorem{theorem}{Theorem}[section]

\newtheorem{lemma}[theorem]{Lemma}
\newtheorem{proposition}[theorem]{Proposition}

\theoremstyle{definition}

\newtheorem{definition}[theorem]{Definition}
\newtheorem{remark}[theorem]{Remark}

\theoremstyle{proposition-definition}
\newtheorem{proposition-definition}[theorem]{Proposition-Definition}

\numberwithin{equation}{section}

\def\ca{{\cal A}} \def\cb{{\cal B}} 
 \def\ce{{\cal E}} \def\cf{{\cal F}}
\def\cg{{\cal G}} \def\ch{{\cal H}} 
  \def\cl{{\cal L}}
  \def\co{{\cal O}}

\def\cs{{\cal S}} \def\ct{{\cal T}} 
\def\cv{{\cal V}}  
 
\def\inj{\hookrightarrow }

\def\gh{GL_n(\HH)\cdot GL_1(\HH)}
\def\gth{GL_n(\HH)\times GL_1(\HH)}
\makeatletter

\def\diagram{\m@th\leftwidth=\z@ \rightwidth=\z@ \topheight=\z@
\botheight=\z@ \setbox\@picbox\hbox\bgroup}

\def\enddiagram{\egroup\wd\@picbox\rightwidth\unitlength
\ht\@picbox\topheight\unitlength \dp\@picbox\botheight\unitlength
\hskip\leftwidth\unitlength\box\@picbox}

\def\bfig{\begin{diagram}}
\def\efig{\end{diagram}}
\newcount\wideness \newcount\leftwidth \newcount\rightwidth
\newcount\highness \newcount\topheight \newcount\botheight

\def\ratchet#1#2{\ifnum#1<#2 \global #1=#2 \fi}

\def\putbox(#1,#2)#3{%
\horsize{\wideness}{#3} \divide\wideness by 2
{\advance\wideness by #1 \ratchet{\rightwidth}{\wideness}}
{\advance\wideness by -#1 \ratchet{\leftwidth}{\wideness}}
\vertsize{\highness}{#3} \divide\highness by 2
{\advance\highness by #2 \ratchet{\topheight}{\highness}}
{\advance\highness by -#2 \ratchet{\botheight}{\highness}}
\put(#1,#2){\makebox(0,0){$#3$}}}

\def\putlbox(#1,#2)#3{%
\horsize{\wideness}{#3}
{\advance\wideness by #1 \ratchet{\rightwidth}{\wideness}}
{\ratchet{\leftwidth}{-#1}}
\vertsize{\highness}{#3} \divide\highness by 2
{\advance\highness by #2 \ratchet{\topheight}{\highness}}
{\advance\highness by -#2 \ratchet{\botheight}{\highness}}
\put(#1,#2){\makebox(0,0)[l]{$#3$}}}

\def\putrbox(#1,#2)#3{%
\horsize{\wideness}{#3}
{\ratchet{\rightwidth}{#1}}
{\advance\wideness by -#1 \ratchet{\leftwidth}{\wideness}}
\vertsize{\highness}{#3} \divide\highness by 2
{\advance\highness by #2 \ratchet{\topheight}{\highness}}
{\advance\highness by -#2 \ratchet{\botheight}{\highness}}
\put(#1,#2){\makebox(0,0)[r]{$#3$}}}

\def\adjust[#1]{} 

\newcount \coefa
\newcount \coefb
\newcount \coefc
\newcount\tempcounta
\newcount\tempcountb
\newcount\tempcountc
\newcount\tempcountd
\newcount\xext
\newcount\yext
\newcount\xoff
\newcount\yoff
\newcount\gap%
\newcount\arrowtypea
\newcount\arrowtypeb
\newcount\arrowtypec
\newcount\arrowtyped
\newcount\arrowtypee
\newcount\height
\newcount\width
\newcount\xpos
\newcount\ypos
\newcount\run
\newcount\rise
\newcount\arrowlength
\newcount\halflength
\newcount\arrowtype
\newdimen\tempdimen
\newdimen\xlen
\newdimen\ylen
\newsavebox{\tempboxa}%
\newsavebox{\tempboxb}%
\newsavebox{\tempboxc}%

\newdimen\w@dth

\def\setw@dth#1#2{\setbox\z@\hbox{\m@th$#1$}\w@dth=\wd\z@
\setbox\@ne\hbox{\m@th$#2$}\ifnum\w@dth<\wd\@ne \w@dth=\wd\@ne \fi
\advance\w@dth by 1.2em}

\def\t@^#1_#2{\allowbreak\def\n@one{#1}\def\n@two{#2}\mathrel
{\setw@dth{#1}{#2}
\mathop{\hbox to \w@dth{\rightarrowfill}}\limits
\ifx\n@one\empty\else ^{\box\z@}\fi
\ifx\n@two\empty\else _{\box\@ne}\fi}}
\def\t@@^#1{\@ifnextchar_{\t@^{#1}}{\t@^{#1}_{}}}
\def\to{\@ifnextchar^{\t@@}{\t@@^{}}}

\def\t@left^#1_#2{\def\n@one{#1}\def\n@two{#2}\mathrel{\setw@dth{#1}{#2}
\mathop{\hbox to \w@dth{\leftarrowfill}}\limits
\ifx\n@one\empty\else ^{\box\z@}\fi
\ifx\n@two\empty\else _{\box\@ne}\fi}}
\def\t@@left^#1{\@ifnextchar_{\t@left^{#1}}{\t@left^{#1}_{}}}
\def\toleft{\@ifnextchar^{\t@@left}{\t@@left^{}}}

\def\two@^#1_#2{\allowbreak
\def\n@one{#1}\def\n@two{#2}\mathrel{\setw@dth{#1}{#2}
\mathop{\vcenter{\lineskip\z@\baselineskip\z@
                 \hbox to \w@dth{\rightarrowfill}%
                 \hbox to \w@dth{\rightarrowfill}}%
       }\limits
\ifx\n@one\empty\else ^{\box\z@}\fi
\ifx\n@two\empty\else _{\box\@ne}\fi}}
\def\tw@@^#1{\@ifnextchar _{\two@^{#1}}{\two@^{#1}_{}}}
\def\two{\@ifnextchar ^{\tw@@}{\tw@@^{}}}

\def\tofr@^#1_#2{\def\n@one{#1}\def\n@two{#2}\mathrel{\setw@dth{#1}{#2}
\mathop{\vcenter{\hbox to \w@dth{\rightarrowfill}\kern-1.7ex
                 \hbox to \w@dth{\leftarrowfill}}%
       }\limits
\ifx\n@one\empty\else ^{\box\z@}\fi
\ifx\n@two\empty\else _{\box\@ne}\fi}}
\def\t@fr@^#1{\@ifnextchar_ {\tofr@^{#1}}{\tofr@^{#1}_{}}}
\def\tofro{\@ifnextchar^ {\t@fr@}{\t@fr@^{}}}

\def\mon{\mathop{\m@th\hbox to
      14.6\P@{\lasyb\char'51\hskip-2.1\P@$\arrext$\hss
$\mathord\rightarrow$}}\limits} 
\def\leftmono{\mathrel{\m@th\hbox to
14.6\P@{$\mathord\leftarrow$\hss$\arrext$\hskip-2.1\P@\lasyb\char'50%
}}\limits} 
\mathchardef\arrext="0200       

\def\settypes(#1,#2,#3){\arrowtypea#1 \arrowtypeb#2 \arrowtypec#3}
\def\settoheight#1#2{\setbox\@tempboxa\hbox{#2}#1\ht\@tempboxa\relax}%
\def\settodepth#1#2{\setbox\@tempboxa\hbox{#2}#1\dp\@tempboxa\relax}%
\def\settokens`#1`#2`#3`#4`{%
     \def\tokena{#1}\def\tokenb{#2}\def\tokenc{#3}\def\tokend{#4}}
\def\setsqparms[#1`#2`#3`#4;#5`#6]{%
\arrowtypea #1
\arrowtypeb #2
\arrowtypec #3
\arrowtyped #4
\width #5
\height #6
}
\def\setpos(#1,#2){\xpos=#1 \ypos#2}

\def\settriparms[#1`#2`#3;#4]{\settripairparms[#1`#2`#3`1`1;#4]}%

\def\settripairparms[#1`#2`#3`#4`#5;#6]{%
\arrowtypea #1
\arrowtypeb #2
\arrowtypec #3
\arrowtyped #4
\arrowtypee #5
\width #6
\height #6
}

\def\resetparms{\settripairparms[1`1`1`1`1;500]\width 500}

\resetparms

\def\mvector(#1,#2)#3{
\put(0,0){\vector(#1,#2){#3}}%
\put(0,0){\vector(#1,#2){26}}%
}
\def\evector(#1,#2)#3{{
\arrowlength #3
\put(0,0){\vector(#1,#2){\arrowlength}}%
\advance \arrowlength by-30
\put(0,0){\vector(#1,#2){\arrowlength}}%
}}

\def\horsize#1#2{%
\settowidth{\tempdimen}{$#2$}%
#1=\tempdimen
\divide #1 by\unitlength
}

\def\vertsize#1#2{%
\settoheight{\tempdimen}{$#2$}%
#1=\tempdimen
\settodepth{\tempdimen}{$#2$}%
\advance #1 by\tempdimen
\divide #1 by\unitlength
}

\def\putvector(#1,#2)(#3,#4)#5#6{{%
\ifnum3<\arrowtype
\putdashvector(#1,#2)(#3,#4)#5\arrowtype
\else
\ifnum\arrowtype<-3
\putdashvector(#1,#2)(#3,#4)#5\arrowtype
\else
\xpos=#1
\ypos=#2
\run=#3
\rise=#4
\arrowlength=#5
\ifnum \arrowtype<0
    \ifnum \run=0
        \advance \ypos by-\arrowlength
    \else
        \tempcounta \arrowlength
        \multiply \tempcounta by\rise
        \divide \tempcounta by\run
        \ifnum\run>0
            \advance \xpos by\arrowlength
            \advance \ypos by\tempcounta
        \else
            \advance \xpos by-\arrowlength
            \advance \ypos by-\tempcounta
        \fi
    \fi
    \multiply \arrowtype by-1
    \multiply \rise by-1
    \multiply \run by-1
\fi
\ifcase \arrowtype
\or \put(\xpos,\ypos){\vector(\run,\rise){\arrowlength}}%
\or \put(\xpos,\ypos){\mvector(\run,\rise)\arrowlength}%
\or \put(\xpos,\ypos){\evector(\run,\rise){\arrowlength}}%
\fi\fi\fi
}}

\def\putsplitvector(#1,#2)#3#4{
\xpos #1
\ypos #2
\arrowtype #4
\halflength #3
\arrowlength #3
\gap 140
\advance \halflength by-\gap
\divide \halflength by2
\ifnum\arrowtype>0
   \ifcase \arrowtype
   \or \put(\xpos,\ypos){\line(0,-1){\halflength}}%
       \advance\ypos by-\halflength
       \advance\ypos by-\gap
       \put(\xpos,\ypos){\vector(0,-1){\halflength}}%
   \or \put(\xpos,\ypos){\line(0,-1)\halflength}%
       \put(\xpos,\ypos){\vector(0,-1)3}%
       \advance\ypos by-\halflength
       \advance\ypos by-\gap
       \put(\xpos,\ypos){\vector(0,-1){\halflength}}%
   \or \put(\xpos,\ypos){\line(0,-1)\halflength}%
       \advance\ypos by-\halflength
       \advance\ypos by-\gap
       \put(\xpos,\ypos){\evector(0,-1){\halflength}}%
   \fi
\else \arrowtype=-\arrowtype
   \ifcase\arrowtype
   \or \advance \ypos by-\arrowlength
       \put(\xpos,\ypos){\line(0,1){\halflength}}%
       \advance\ypos by\halflength
       \advance\ypos by\gap
       \put(\xpos,\ypos){\vector(0,1){\halflength}}%
   \or \advance \ypos by-\arrowlength
       \put(\xpos,\ypos){\line(0,1)\halflength}%
       \put(\xpos,\ypos){\vector(0,1)3}%
       \advance\ypos by\halflength
       \advance\ypos by\gap
       \put(\xpos,\ypos){\vector(0,1){\halflength}}%
   \or \advance \ypos by-\arrowlength
       \put(\xpos,\ypos){\line(0,1)\halflength}%
       \advance\ypos by\halflength
       \advance\ypos by\gap
       \put(\xpos,\ypos){\evector(0,1){\halflength}}%
   \fi
\fi
}

\def\putmorphism(#1)(#2,#3)[#4`#5`#6]#7#8#9{{%
\run #2
\rise #3
\ifnum\rise=0
  \puthmorphism(#1)[#4`#5`#6]{#7}{#8}#9%
\else\ifnum\run=0
  \putvmorphism(#1)[#4`#5`#6]{#7}{#8}#9%
\else
\setpos(#1)%
\arrowlength #7
\arrowtype #8
\ifnum\run=0
\else\ifnum\rise=0
\else
\ifnum\run>0
    \coefa=1
\else
   \coefa=-1
\fi
\ifnum\arrowtype>0
   \coefb=0
   \coefc=-1
\else
   \coefb=\coefa
   \coefc=1
   \arrowtype=-\arrowtype
\fi
\width=2
\multiply \width by\run
\divide \width by\rise
\ifnum \width<0  \width=-\width\fi
\advance\width by60
\if l#9 \width=-\width\fi
\putbox(\xpos,\ypos){#4}
{\multiply \coefa by\arrowlength
\advance\xpos by\coefa
\multiply \coefa by\rise
\divide \coefa by\run
\advance \ypos by\coefa
\putbox(\xpos,\ypos){#5} }%
{\multiply \coefa by\arrowlength
\divide \coefa by2
\advance \xpos by\coefa
\advance \xpos by\width
\multiply \coefa by\rise
\divide \coefa by\run
\advance \ypos by\coefa
\if l#9%
   \putrbox(\xpos,\ypos){#6}%
\else\if r#9%
   \putlbox(\xpos,\ypos){#6}%
\fi\fi }%
{\multiply \rise by-\coefc
\multiply \run by-\coefc
\multiply \coefb by\arrowlength
\advance \xpos by\coefb
\multiply \coefb by\rise
\divide \coefb by\run
\advance \ypos by\coefb
\multiply \coefc by70
\advance \ypos by\coefc
\multiply \coefc by\run
\divide \coefc by\rise
\advance \xpos by\coefc
\multiply \coefa by140
\multiply \coefa by\run
\divide \coefa by\rise
\advance \arrowlength by\coefa
\ifcase\arrowtype
\or \put(\xpos,\ypos){\vector(\run,\rise){\arrowlength}}%
\or \put(\xpos,\ypos){\mvector(\run,\rise){\arrowlength}}%
\or \put(\xpos,\ypos){\evector(\run,\rise){\arrowlength}}%
\fi}\fi\fi\fi\fi}}

\newcount\numbdashes \newcount\lengthdash \newcount\increment

\def\howmanydashes{
\numbdashes=\arrowlength \lengthdash=40
\divide\numbdashes by \lengthdash
\lengthdash=\arrowlength
\divide\lengthdash by \numbdashes
\increment=\lengthdash
\multiply\lengthdash by 3
\divide\lengthdash by 5
}

\def\putdashvector(#1)(#2,#3)#4#5{%
\ifnum#3=0 \putdashhvector(#1){#4}#5
\else
\ifnum#2=0
\putdashvvector(#1){#4}#5\fi\fi}

\def\putdashhvector(#1,#2)#3#4{{%
\arrowlength=#3 \howmanydashes
\multiput(#1,#2)(\increment,0){\numbdashes}%
{\vrule height .4pt width \lengthdash\unitlength}
\arrowtype=#4 \xpos=#1
\ifnum\arrowtype<0 \advance\arrowtype by 7 \fi
\ifcase\arrowtype
\or \advance\xpos by 10
    \put(\xpos,#2){\vector(-1,0){\lengthdash}}
    \advance\xpos by 40
    \put(\xpos,#2){\vector(-1,0){\lengthdash}}
\or \advance \xpos by 10
    \put(\xpos,#2){\vector(-1,0){\lengthdash}}
    \advance\xpos by  \arrowlength
    \advance\xpos by  -50
    \put(\xpos,#2){\vector(-1,0){\lengthdash}}
\or \advance\xpos by 10
    \put(\xpos,#2){\vector(-1,0){\lengthdash}}
\or \advance\xpos by \arrowlength
    \advance\xpos by -\lengthdash
    \put(\xpos,#2){\vector(1,0){\lengthdash}}
\or {\advance\xpos by 10
    \put(\xpos,#2){\vector(1,0){\lengthdash}}}
    \advance\xpos by \arrowlength
    \advance\xpos by -\lengthdash
    \put(\xpos,#2){\vector(1,0){\lengthdash}}
\or \advance\xpos by \arrowlength
    \advance\xpos by -\lengthdash
    \put(\xpos,#2){\vector(1,0){\lengthdash}}
    \advance\xpos by -40
    \put(\xpos,#2){\vector(1,0){\lengthdash}}
   \fi
}}

\def\putdashvvector(#1,#2)#3#4{{%
\arrowlength=#3 \howmanydashes
\ypos=#2 \advance\ypos by -\arrowlength
\multiput(#1,#2)(0,\increment){\numbdashes}%
    {\vrule width .4pt height \lengthdash\unitlength}
\arrowtype=#4 \ypos=#2
\ifnum\arrowtype<0 \advance\arrowtype by 7 \fi
\ifcase\arrowtype
\or \advance\ypos by \arrowlength \advance\ypos by -40
    \put(#1,\ypos){\vector(0,1){\lengthdash}}
    \advance\ypos by -40
    \put(#1,\ypos){\vector(0,1){\lengthdash}}
\or \advance\ypos by 10
    \put(#1,\ypos){\vector(0,1){\lengthdash}}
    \advance\ypos by \arrowlength \advance\ypos by -40
    \put(#1,\ypos){\vector(0,1){\lengthdash}}
\or \advance\ypos by \arrowlength \advance\ypos by -40
    \put(#1,\ypos){\vector(0,1){\lengthdash}}
\or \advance\ypos by 10
    \put(#1,\ypos){\vector(0,-1){\lengthdash}}
\or \advance\ypos by 10
    \put(#1,\ypos){\vector(0,-1){\lengthdash}}
    \advance\ypos by \arrowlength \advance\ypos by -40
    \put(#1,\ypos){\vector(0,-1){\lengthdash}}
\or \advance\ypos by 10
    \put(#1,\ypos){\vector(0,-1){\lengthdash}}
    \advance\ypos by 40
    \put(#1,\ypos){\vector(0,-1){\lengthdash}}
\fi
}}

\def\puthmorphism(#1,#2)[#3`#4`#5]#6#7#8{{%
\xpos #1
\ypos #2
\width #6
\arrowlength #6
\arrowtype=#7
\putbox(\xpos,\ypos){#3\vphantom{#4}}%
{\advance \xpos by\arrowlength
\putbox(\xpos,\ypos){\vphantom{#3}#4}}%
\horsize{\tempcounta}{#3}%
\horsize{\tempcountb}{#4}%
\divide \tempcounta by2
\divide \tempcountb by2
\advance \tempcounta by30
\advance \tempcountb by30
\advance \xpos by\tempcounta
\advance \arrowlength by-\tempcounta
\advance \arrowlength by-\tempcountb
\putvector(\xpos,\ypos)(1,0)\arrowlength\arrowtype
\divide \arrowlength by2
\advance \xpos by\arrowlength
\vertsize{\tempcounta}{#5}%
\divide\tempcounta by2
\advance \tempcounta by20
\if a#8 %
   \advance \ypos by\tempcounta
   \putbox(\xpos,\ypos){#5}%
\else
   \advance \ypos by-\tempcounta
   \putbox(\xpos,\ypos){#5}%
\fi}}

\def\putvmorphism(#1,#2)[#3`#4`#5]#6#7#8{{%
\xpos #1
\ypos #2
\arrowlength #6
\arrowtype #7
\settowidth{\xlen}{$#5$}%
\putbox(\xpos,\ypos){#3}%
{\advance \ypos by-\arrowlength
\putbox(\xpos,\ypos){#4}}%
{\advance\arrowlength by-140
\advance \ypos by-70
\ifdim\xlen>0pt
   \if m#8%
      \putsplitvector(\xpos,\ypos)\arrowlength\arrowtype
   \else
   \putvector(\xpos,\ypos)(0,-1)\arrowlength\arrowtype
   \fi
\else
   \putvector(\xpos,\ypos)(0,-1)\arrowlength\arrowtype
\fi}%
\ifdim\xlen>0pt
   \divide \arrowlength by2
   \advance\ypos by-\arrowlength
   \if l#8%
      \advance \xpos by-40
      \putrbox(\xpos,\ypos){#5}%
   \else\if r#8%
      \advance \xpos by40
      \putlbox(\xpos,\ypos){#5}%
   \else
      \putbox(\xpos,\ypos){#5}%
   \fi\fi
\fi
}}

\def\putsquarep<#1>(#2)[#3;#4`#5`#6`#7]{{%
\setsqparms[#1]%
\setpos(#2)%
\settokens`#3`%
\puthmorphism(\xpos,\ypos)[\tokenc`\tokend`{#7}]{\width}{\arrowtyped}b%
\advance\ypos by \height
\puthmorphism(\xpos,\ypos)[\tokena`\tokenb`{#4}]{\width}{\arrowtypea}a%
\putvmorphism(\xpos,\ypos)[``{#5}]{\height}{\arrowtypeb}l%
\advance\xpos by \width
\putvmorphism(\xpos,\ypos)[``{#6}]{\height}{\arrowtypec}r%
}}

\def\putsquare{\@ifnextchar <{\putsquarep}{\putsquarep%
   <\arrowtypea`\arrowtypeb`\arrowtypec`\arrowtyped;\width`\height>}}
\def\square{\@ifnextchar< {\squarep}{\squarep
   <\arrowtypea`\arrowtypeb`\arrowtypec`\arrowtyped;\width`\height>}}
\def\squarep<#1>[#2`#3`#4`#5;#6`#7`#8`#9]{{
\setsqparms[#1]
\diagram
\putsquarep<\arrowtypea`\arrowtypeb`\arrowtypec`
\arrowtyped;\width`\height>
(0,0)[#2`#3`#4`{#5};#6`#7`#8`{#9}]
\enddiagram
}}                                                 
\def\putptrianglep<#1>(#2,#3)[#4`#5`#6;#7`#8`#9]{{%
\settriparms[#1]%
\xpos=#2 \ypos=#3
\advance\ypos by \height
\puthmorphism(\xpos,\ypos)[#4`#5`{#7}]{\height}{\arrowtypea}a%
\putvmorphism(\xpos,\ypos)[`#6`{#8}]{\height}{\arrowtypeb}l%
\advance\xpos by\height
\putmorphism(\xpos,\ypos)(-1,-1)[``{#9}]{\height}{\arrowtypec}r%
}}

\def\putptriangle{\@ifnextchar <{\putptrianglep}{\putptrianglep
   <\arrowtypea`\arrowtypeb`\arrowtypec;\height>}}
\def\ptriangle{\@ifnextchar <{\ptrianglep}{\ptrianglep
   <\arrowtypea`\arrowtypeb`\arrowtypec;\height>}}
\def\ptrianglep<#1>[#2`#3`#4;#5`#6`#7]{{
\settriparms[#1]
\diagram
\putptrianglep<\arrowtypea`\arrowtypeb`
\arrowtypec;\height>
(0,0)[#2`#3`#4;#5`#6`{#7}]
\enddiagram
}}                                            

\def\putqtrianglep<#1>(#2,#3)[#4`#5`#6;#7`#8`#9]{{%
\settriparms[#1]%
\xpos=#2 \ypos=#3
\advance\ypos by\height
\puthmorphism(\xpos,\ypos)[#4`#5`{#7}]{\height}{\arrowtypea}a%
\putmorphism(\xpos,\ypos)(1,-1)[``{#8}]{\height}{\arrowtypeb}l%
\advance\xpos by\height
\putvmorphism(\xpos,\ypos)[`#6`{#9}]{\height}{\arrowtypec}r%
}}

\def\putqtriangle{\@ifnextchar <{\putqtrianglep}{\putqtrianglep
   <\arrowtypea`\arrowtypeb`\arrowtypec;\height>}}
\def\qtriangle{\@ifnextchar <{\qtrianglep}{\qtrianglep
   <\arrowtypea`\arrowtypeb`\arrowtypec;\height>}}
\def\qtrianglep<#1>[#2`#3`#4;#5`#6`#7]{{
\settriparms[#1]
\width=\height                                
\diagram
\putqtrianglep<\arrowtypea`\arrowtypeb`
\arrowtypec;\height>
(0,0)[#2`#3`#4;#5`#6`{#7}]
\enddiagram
}}

\def\putdtrianglep<#1>(#2,#3)[#4`#5`#6;#7`#8`#9]{{%
\settriparms[#1]%
\xpos=#2 \ypos=#3
\puthmorphism(\xpos,\ypos)[#5`#6`{#9}]{\height}{\arrowtypec}b%
\advance\xpos by \height \advance\ypos by\height
\putmorphism(\xpos,\ypos)(-1,-1)[``{#7}]{\height}{\arrowtypea}l%
\putvmorphism(\xpos,\ypos)[#4``{#8}]{\height}{\arrowtypeb}r%
}}

\def\putdtriangle{\@ifnextchar <{\putdtrianglep}{\putdtrianglep
   <\arrowtypea`\arrowtypeb`\arrowtypec;\height>}}
\def\dtriangle{\@ifnextchar <{\dtrianglep}{\dtrianglep
   <\arrowtypea`\arrowtypeb`\arrowtypec;\height>}}
\def\dtrianglep<#1>[#2`#3`#4;#5`#6`#7]{{
\settriparms[#1]
\width=\height                                
\diagram
\putdtrianglep<\arrowtypea`\arrowtypeb`
\arrowtypec;\height>
(0,0)[#2`#3`#4;#5`#6`{#7}]
\enddiagram
}}

\def\putbtrianglep<#1>(#2,#3)[#4`#5`#6;#7`#8`#9]{{%
\settriparms[#1]%
\xpos=#2 \ypos=#3
\puthmorphism(\xpos,\ypos)[#5`#6`{#9}]{\height}{\arrowtypec}b%
\advance\ypos by\height
\putmorphism(\xpos,\ypos)(1,-1)[``{#8}]{\height}{\arrowtypeb}r%
\putvmorphism(\xpos,\ypos)[#4``{#7}]{\height}{\arrowtypea}l%
}}

\def\putbtriangle{\@ifnextchar <{\putbtrianglep}{\putbtrianglep
   <\arrowtypea`\arrowtypeb`\arrowtypec;\height>}}
\def\btriangle{\@ifnextchar <{\btrianglep}{\btrianglep
   <\arrowtypea`\arrowtypeb`\arrowtypec;\height>}}
\def\btrianglep<#1>[#2`#3`#4;#5`#6`#7]{{
\settriparms[#1]
\width=\height                               
\diagram
\putbtrianglep<\arrowtypea`\arrowtypeb`
\arrowtypec;\height>
(0,0)[#2`#3`#4;#5`#6`{#7}]
\enddiagram
}}

\def\putAtrianglep<#1>(#2,#3)[#4`#5`#6;#7`#8`#9]{{%
\settriparms[#1]%
\xpos=#2 \ypos=#3
{\multiply \height by2
\puthmorphism(\xpos,\ypos)[#5`#6`{#9}]{\height}{\arrowtypec}b}%
\advance\xpos by\height \advance\ypos by\height
\putmorphism(\xpos,\ypos)(-1,-1)[#4``{#7}]{\height}{\arrowtypea}l%
\putmorphism(\xpos,\ypos)(1,-1)[``{#8}]{\height}{\arrowtypeb}r%
}}

\def\putAtriangle{\@ifnextchar <{\putAtrianglep}{\putAtrianglep
   <\arrowtypea`\arrowtypeb`\arrowtypec;\height>}}
\def\Atriangle{\@ifnextchar <{\Atrianglep}{\Atrianglep
   <\arrowtypea`\arrowtypeb`\arrowtypec;\height>}}
\def\Atrianglep<#1>[#2`#3`#4;#5`#6`#7]{{
\settriparms[#1]
\width=\height                                     
\diagram
\putAtrianglep<\arrowtypea`\arrowtypeb`
\arrowtypec;\height>
(0,0)[#2`#3`#4;#5`#6`{#7}]
\enddiagram
}}

\def\putAtrianglepairp<#1>(#2)[#3;#4`#5`#6`#7`#8]{{%
\settripairparms[#1]%
\setpos(#2)%
\settokens`#3`%
\puthmorphism(\xpos,\ypos)[\tokenb`\tokenc`{#7}]{\height}{\arrowtyped}b%
\advance\xpos by\height
\puthmorphism(\xpos,\ypos)[\phantom{\tokenc}`\tokend`{#8}]%
{\height}{\arrowtypee}b%
\advance\ypos by\height
\putmorphism(\xpos,\ypos)(-1,-1)[\tokena``{#4}]{\height}{\arrowtypea}l%
\putvmorphism(\xpos,\ypos)[``{#5}]{\height}{\arrowtypeb}m%
\putmorphism(\xpos,\ypos)(1,-1)[``{#6}]{\height}{\arrowtypec}r%
}}

\def\putAtrianglepair{\@ifnextchar <{\putAtrianglepairp}{\putAtrianglepairp%
   <\arrowtypea`\arrowtypeb`\arrowtypec`\arrowtyped`\arrowtypee;\height>}}
\def\Atrianglepair{\@ifnextchar <{\Atrianglepairp}{\Atrianglepairp%
   <\arrowtypea`\arrowtypeb`\arrowtypec`\arrowtyped`\arrowtypee;\height>}}

\def\Atrianglepairp<#1>[#2;#3`#4`#5`#6`#7]{{
\settripairparms[#1]
\settokens`#2`
\width=\height                                
\diagram
\putAtrianglepairp                            
<\arrowtypea`\arrowtypeb`\arrowtypec`
\arrowtyped`\arrowtypee;\height>
(0,0)[{#2};#3`#4`#5`#6`{#7}]
\enddiagram
}}

\def\putVtrianglep<#1>(#2,#3)[#4`#5`#6;#7`#8`#9]{{%
\settriparms[#1]%
\xpos=#2 \ypos=#3
\advance\ypos by\height
{\multiply\height by2
\puthmorphism(\xpos,\ypos)[#4`#5`{#7}]{\height}{\arrowtypea}a}%
\putmorphism(\xpos,\ypos)(1,-1)[`#6`{#8}]{\height}{\arrowtypeb}l%
\advance\xpos by\height
\advance\xpos by\height
\putmorphism(\xpos,\ypos)(-1,-1)[``{#9}]{\height}{\arrowtypec}r%
}}

\def\putVtriangle{\@ifnextchar <{\putVtrianglep}{\putVtrianglep
   <\arrowtypea`\arrowtypeb`\arrowtypec;\height>}}
\def\Vtriangle{\@ifnextchar <{\Vtrianglep}{\Vtrianglep
   <\arrowtypea`\arrowtypeb`\arrowtypec;\height>}}
\def\Vtrianglep<#1>[#2`#3`#4;#5`#6`#7]{{
\settriparms[#1]
\width=\height                                 
\diagram
\putVtrianglep<\arrowtypea`\arrowtypeb`
\arrowtypec;\height>
(0,0)[#2`#3`#4;#5`#6`{#7}]
\enddiagram
}}

\def\putVtrianglepairp<#1>(#2)[#3;#4`#5`#6`#7`#8]{{
\settripairparms[#1]%
\setpos(#2)%
\settokens`#3`%
\advance\ypos by\height
\putmorphism(\xpos,\ypos)(1,-1)[`\tokend`{#6}]{\height}{\arrowtypec}l%
\puthmorphism(\xpos,\ypos)[\tokena`\tokenb`{#4}]{\height}{\arrowtypea}a%
\advance\xpos by\height
\puthmorphism(\xpos,\ypos)[\phantom{\tokenb}`\tokenc`{#5}]%
{\height}{\arrowtypeb}a%
\putvmorphism(\xpos,\ypos)[``{#7}]{\height}{\arrowtyped}m%
\advance\xpos by\height
\putmorphism(\xpos,\ypos)(-1,-1)[``{#8}]{\height}{\arrowtypee}r%
}}

\def\putVtrianglepair{\@ifnextchar <{\putVtrianglepairp}{\putVtrianglepairp%
    <\arrowtypea`\arrowtypeb`\arrowtypec`\arrowtyped`\arrowtypee;\height>}}
\def\Vtrianglepair{\@ifnextchar <{\Vtrianglepairp}{\Vtrianglepairp%
    <\arrowtypea`\arrowtypeb`\arrowtypec`\arrowtyped`\arrowtypee;\height>}}
\def\Vtrianglepairp<#1>[#2;#3`#4`#5`#6`#7]{{
\settripairparms[#1]
\settokens`#2`
\diagram
\putVtrianglepairp                             
<\arrowtypea`\arrowtypeb`\arrowtypec`
\arrowtyped`\arrowtypee;\height>
(0,0)[{#2};#3`#4`#5`#6`{#7}]
\enddiagram
}}

\def\putCtrianglep<#1>(#2,#3)[#4`#5`#6;#7`#8`#9]{{%
\settriparms[#1]%
\xpos=#2 \ypos=#3
\advance\ypos by\height
\putmorphism(\xpos,\ypos)(1,-1)[``{#9}]{\height}{\arrowtypec}l%
\advance\xpos by\height
\advance\ypos by\height
\putmorphism(\xpos,\ypos)(-1,-1)[#4`#5`{#7}]{\height}{\arrowtypea}l%
{\multiply\height by 2
\putvmorphism(\xpos,\ypos)[`#6`{#8}]{\height}{\arrowtypeb}r}%
}}

\def\putCtriangle{\@ifnextchar <{\putCtrianglep}{\putCtrianglep
    <\arrowtypea`\arrowtypeb`\arrowtypec;\height>}}
\def\Ctriangle{\@ifnextchar <{\Ctrianglep}{\Ctrianglep
    <\arrowtypea`\arrowtypeb`\arrowtypec;\height>}}
\def\Ctrianglep<#1>[#2`#3`#4;#5`#6`#7]{{
\settriparms[#1]
\width=\height                               
\diagram
\putCtrianglep<\arrowtypea`\arrowtypeb`
\arrowtypec;\height>
(0,0)[#2`#3`#4;#5`#6`{#7}]
\enddiagram
}}                                           
\def\putDtrianglep<#1>(#2,#3)[#4`#5`#6;#7`#8`#9]{{%
\settriparms[#1]%
\xpos=#2 \ypos=#3
\advance\xpos by\height \advance\ypos by\height
\putmorphism(\xpos,\ypos)(-1,-1)[``{#9}]{\height}{\arrowtypec}r%
\advance\xpos by-\height \advance\ypos by\height
\putmorphism(\xpos,\ypos)(1,-1)[`#5`{#8}]{\height}{\arrowtypeb}r%
{\multiply\height by 2
\putvmorphism(\xpos,\ypos)[#4`#6`{#7}]{\height}{\arrowtypea}l}%
}}

\def\putDtriangle{\@ifnextchar <{\putDtrianglep}{\putDtrianglep
    <\arrowtypea`\arrowtypeb`\arrowtypec;\height>}}
\def\Dtriangle{\@ifnextchar <{\Dtrianglep}{\Dtrianglep
   <\arrowtypea`\arrowtypeb`\arrowtypec;\height>}}
\def\Dtrianglep<#1>[#2`#3`#4;#5`#6`#7]{{
\settriparms[#1]
\width=\height                              
\diagram
\putDtrianglep<\arrowtypea`\arrowtypeb`
\arrowtypec;\height>
(0,0)[#2`#3`#4;#5`#6`{#7}]
\enddiagram
}}                                          
\def\setrecparms[#1`#2]{\width=#1 \height=#2}%

\def\recursep<#1`#2>[#3;#4`#5`#6`#7`#8]{{\m@th
\width=#1 \height=#2
\settokens`#3`
\settowidth{\tempdimen}{$\tokena$}
\ifdim\tempdimen=0pt
  \savebox{\tempboxa}{\hbox{$\tokenb$}}%
  \savebox{\tempboxb}{\hbox{$\tokend$}}%
  \savebox{\tempboxc}{\hbox{$#6$}}%
\else
  \savebox{\tempboxa}{\hbox{$\hbox{$\tokena$}\times\hbox{$\tokenb$}$}}%
  \savebox{\tempboxb}{\hbox{$\hbox{$\tokena$}\times\hbox{$\tokend$}$}}%
  \savebox{\tempboxc}{\hbox{$\hbox{$\tokena$}\times\hbox{$#6$}$}}%
\fi
\ypos=\height
\divide\ypos by 2
\xpos=\ypos
\advance\xpos by \width
\bfig
\putCtrianglep<-1`1`1;\ypos>(0,0)[`\tokenc`;#5`#6`{#7}]%
\puthmorphism(\ypos,0)[\tokend`\usebox{\tempboxb}`{#8}]{\width}{-1}b%
\puthmorphism(\ypos,\height)[\tokenb`\usebox{\tempboxa}`{#4}]{\width}{-1}a%
\advance\ypos by \width
\putvmorphism(\ypos,\height)[``\usebox{\tempboxc}]{\height}1r%
\efig
}}

\def\recurse{\@ifnextchar <{\recursep}{\recursep<\width`\height>}}

\def\puttwohmorphisms(#1,#2)[#3`#4;#5`#6]#7#8#9{{%
\puthmorphism(#1,#2)[#3`#4`]{#7}0a
\ypos=#2
\advance\ypos by 20
\puthmorphism(#1,\ypos)[\phantom{#3}`\phantom{#4}`#5]{#7}{#8}a
\advance\ypos by -40
\puthmorphism(#1,\ypos)[\phantom{#3}`\phantom{#4}`#6]{#7}{#9}b
}}

\def\puttwovmorphisms(#1,#2)[#3`#4;#5`#6]#7#8#9{{%
\putvmorphism(#1,#2)[#3`#4`]{#7}0a
\xpos=#1
\advance\xpos by -20
\putvmorphism(\xpos,#2)[\phantom{#3}`\phantom{#4}`#5]{#7}{#8}l
\advance\xpos by 40
\putvmorphism(\xpos,#2)[\phantom{#3}`\phantom{#4}`#6]{#7}{#9}r
}}

\def\puthcoequalizer(#1)[#2`#3`#4;#5`#6`#7]#8#9{{%
\setpos(#1)%
\puttwohmorphisms(\xpos,\ypos)[#2`#3;#5`#6]{#8}11%
\advance\xpos by #8
\puthmorphism(\xpos,\ypos)[\phantom{#3}`#4`#7]{#8}1{#9}
}}

\def\putvcoequalizer(#1)[#2`#3`#4;#5`#6`#7]#8#9{{%
\setpos(#1)%
\puttwovmorphisms(\xpos,\ypos)[#2`#3;#5`#6]{#8}11%
\advance\ypos by -#8
\putvmorphism(\xpos,\ypos)[\phantom{#3}`#4`#7]{#8}1{#9}
}}

\def\putthreehmorphisms(#1)[#2`#3;#4`#5`#6]#7(#8)#9{{%
\setpos(#1) \settypes(#8)
\if a#9 %
     \vertsize{\tempcounta}{#5}%
     \vertsize{\tempcountb}{#6}%
     \ifnum \tempcounta<\tempcountb \tempcounta=\tempcountb \fi
\else
     \vertsize{\tempcounta}{#4}%
     \vertsize{\tempcountb}{#5}%
     \ifnum \tempcounta<\tempcountb \tempcounta=\tempcountb \fi
\fi
\advance \tempcounta by 60
\puthmorphism(\xpos,\ypos)[#2`#3`#5]{#7}{\arrowtypeb}{#9}
\advance\ypos by \tempcounta
\puthmorphism(\xpos,\ypos)[\phantom{#2}`\phantom{#3}`#4]{#7}{\arrowtypea}{#9}
\advance\ypos by -\tempcounta \advance\ypos by -\tempcounta
\puthmorphism(\xpos,\ypos)[\phantom{#2}`\phantom{#3}`#6]{#7}{\arrowtypec}{#9}
}}

\def\setarrowtoks[#1`#2`#3`#4`#5`#6]{%
\def\toka{#1}
\def\tokb{#2}
\def\tokc{#3}
\def\tokd{#4}
\def\toke{#5}
\def\tokf{#6}
}
\def\hex{\@ifnextchar <{\hexp}{\hexp<1000`400>}}
\def\hexp<#1`#2>[#3`#4`#5`#6`#7`#8;#9]{%
\setarrowtoks[#9]
\yext=#2 \advance \yext by #2
\xext=#1 \advance\xext by \yext
\bfig
\putCtriangle<-1`0`1;#2>(0,0)[`#5`;\tokb``\tokd]
\xext=#1 \yext=#2 \advance \yext by #2
\putsquare<1`0`0`1;\xext`\yext>(#2,0)[#3`#4`#7`#8;\toka```\tokf]
\advance \xext by #2
\putDtriangle<0`1`-1;#2>(\xext,0)[`#6`;`\tokc`\toke]
\efig
}
\makeatother

\begin{document}
\maketitle
\begin{abstract}
On any quaternionic manifold of dimension greater than 4 a class of
plurisubharmonic functions (or, rather, sections of an appropriate
line bundle) is introduced. Then a Monge-Amp\`ere operator is
defined. It is shown that it satisfies a version of theorems of A.
D. Alexandrov and Chern-Levine-Nirenberg. For more special classes
of manifolds analogous results were previously obtained in
\cite{alesker-bsm-03} for the flat quaternionic space $\HH^n$ and in
\cite{alesker-verbitsky-06} for hypercomplex manifolds. One of the
new technical aspects of the present paper is the systematic use of
the Baston differential operators, for which we also prove a new
multiplicativity property.
\end{abstract}

\tableofcontents \setcounter{section}{-1}

\section{Introduction.}\label{S:introduction}
In the recent years the classical theory of plurisubharmonic
functions of complex variables has been generalized in several
directions. These generalizations are in some respects analogous to
the complex case and the real one  (i.e. the theory of convex
functions), but nevertheless reflect rather different geometry
behind.

The author \cite{alesker-bsm-03} and independently at the same time
G. Henkin \cite{henkin} have introduced and studied a class of
plurisubharmonic functions of quaternionic variables on the flat
quaternionic space $\HH^n$. This class was studied further in
\cite{alesker-adv-05} where the author has also obtained
applications to the theory of valuations on convex sets. Analogous,
though geometrically different, results with applications to the
valuations theory were obtained in \cite{alesker-octonions} for the
case of plurisubharmonic functions of two octonionic variables.

A class of plurisubharmonic functions on hypercomplex manifolds was
introduced by M. Verbitsky and the author
\cite{alesker-verbitsky-06}; in the special case of the flat
hypercomplex manifold $\HH^n$ this class coincides with the above
mentioned one. Also in \cite{alesker-verbitsky-06} an interesting
geometric interpretation of strictly plurisubharmonic functions on
hypercomplex manifolds was obtained: locally they are precisely the
potentials for a special class of Riemannian metrics called
Hyper-K\"ahler with Torsion (HKT). This is obviously analogous to
the well known interpretation of strictly plurisubharmonic functions
of complex variables as local potentials of K\"ahler metrics.

In the above mentioned papers an important role was played by
quaternionic (and octonionic in \cite{alesker-octonions}) versions
of the Hessian and the Monge-Amp\`ere operator. They were applied
further for the theory of quaternionic Monge-Amp\`ere equations in
\cite{alesker-jga-03} in the flat case and in
\cite{alesker-verbitsky-09}, \cite{verbitsky-ma} on hypercomplex
manifolds in the context of HKT-geometry.

\hfill

In the recent series of articles
\cite{harvey-lawson1}-\cite{harvey-lawson3} Harvey and Lawson have
developed another approach to pluripotential theory in the context
of calibrated geometries. In various special cases it partly
overlaps with the above mentioned approach. For example on the flat
space $\HH^n$ the two approaches lead to the same class of
plurisubharmonic functions.

\hfill

In this paper we introduce a class of plurisubharmonic functions (or
more precisely, sections of certain specific line bundle) on an
arbitrary quaternionic manifold. Quaternionic manifolds were
introduced independently by S. Salamon \cite{salamon-82} and L.
B\'erard-Bergery (see Ch. 14 in the book \cite{besse}). They carry
quite rich structures, for instance admit the twistor space
\cite{salamon-82} and various canonically associated differential
operators, e.g. the Baston operators (see \cite{baston-92} and
Section \ref{S:penrose-baston} below). Important special cases of
quaternionic manifolds are the flat space $\HH^n$, the quaternionic
projective space $\HH\PP^n$, hypercomplex manifolds (in particular,
hyper-K\"ahler manifolds), and quaternionic K\"ahler manifolds. For
more examples of quaternionic manifolds we refer to Ch. 14 of the
book \cite{besse} and to \cite{joyce}.

We introduce a Monge-Amp\`ere operator on quaternionic manifolds and
prove a version of theorems of A.D. Aleksandrov \cite{alexandrov-58}
and Chern-Levine-Nirenberg \cite{chern-levine-nirenberg}. In the
special cases of the flat space and hypercomplex manifolds analogous
results were obtained in \cite{alesker-bsm-03} and
\cite{alesker-verbitsky-06} respectively. Formally speaking, the
theory in the flat case in \cite{alesker-bsm-03} is indeed a special
case of the theory developed in this paper (see Section
\ref{S:hypercomplex}), as well as it is a special case of the theory
\cite{alesker-verbitsky-06} in the hypercomplex case. However for a
general hypercomplex manifold the theory of
\cite{alesker-verbitsky-06} is not a special case of the theory of
this paper, at least in the case when the holonomy of the Obata
connection is not contained in $SL_n(\HH)$: for example
plurisubharmonic sections belong to different line bundles in the
two theories. It would be of interest to make a more detailed
comparison of the two approaches, in particular when the holonomy of
the Obata connection of a hypercomplex manifold is contained in
$SL_n(\HH)$.

\hfill

The new step of the present paper is the use of the Baston
differential operators; one of them turns out to be the right
quaternionic analogue of the Hessian on a general quaternionic
manifold.

\begin{remark}
One new special case covered by this paper in comparison to
\cite{alesker-bsm-03}, \cite{alesker-verbitsky-06} is the case of
quaternionic K\"ahler manifolds. Since such manifolds can be
considered from the point of view of calibrated geometry, they are
also covered by the Harvey-Lawson theory. It turns out that in
general the Harvey-Lawson class of plurisubharmonic functions is
different from that introduced here, e.g. in the case of
quaternionic projective space $\HH\PP^n$ equipped with the
Fubini-Study metric. Nevertheless when the metric is flat, e.g. on
$\HH^n$, the two classes do coincide.
\end{remark}

\hfill

Let us describe the main results in greater detail. The first step
is to introduce the right version of of quaternionic Hessian on a
general quaternionic manifold $M^{4n}$. We assume throughout the
article that $n>1$; the case $n=1$ is more elementary but somewhat
exceptional for quaternionic manifolds. We claim that this is a
differential operator $\Delta$ of second order which was introduced
for general quaternionic manifolds by Baston \cite{baston-92} for
completely different reasons. The operator $\Delta$ is defined on
smooth sections of a real vector bundle over $M$ which is denoted by
$(\det\ch^*_0)_\RR$ (for the moment, take it as a single notation),
and takes values in a vector bundle denoted by
$\wedge^2\ce^*_0[-2]_\RR$. It is analogous to the operator $dd^c$ on
a complex manifold. The Baston's construction is discussed in
Section \ref{S:penrose-baston}. It is rather involved and uses the
twistor space of $M$ and the Penrose transform. Notice that in the
flat case essentially the same construction was invented much
earlier by Gindikin and Henkin \cite{gindikin-henkin}. It is shown
in Section \ref{S:hypercomplex} that on the flat space $\HH^n$ the
operator $\Delta$ coincides with the quaternionic Hessian introduced
in \cite{alesker-bsm-03} whose construction was elementary.

Next we define the notion of positivity in the fibers of the bundle
$\wedge^2\ce^*_0[-2]_\RR$ (Section \ref{S:positive-currents}). Then
we call a $C^2$-smooth section $h$ of $(\det\ch_0^*)_\RR$ to be
plurisubharmonic if $\Delta h$ is non-negative. \footnote{In fact we
define plurisubharmonicity in a slightly greater generality of just
continuous sections.}

Now let us state the main theorem of the paper. For the sake of
simplicity of exposition we do it here in a weaker form; we refer to
Section \ref{S:psh} for a complete statement. We have the natural
wedge product
$$\wedge^p\ce^*_0[i]_\RR\otimes \wedge^q\ce_0^*[j]_\RR\to
\wedge^{p+q}\ce_0^*[i+j]_\RR.$$ Using this product we define the
Monge-Amp\`ere operator on sections of $(\det\ch^*_0)_\RR$ by
$$h\mapsto (\Delta h)^n$$
where $4n$ is the real dimension of the manifold $M$ as previously.
Thus the Monge-Amp\`ere operator takes values in the real line
bundle $\det\ce_0^*[-2n]_\RR$. Then we show that if $\{h_N\}$ is a
sequence of $C^2$-smooth plurisubharmonic sections which converges
to a $C^2$-section $h$ in the $C^0$-topology (not necessarily in the
$C^2$-topology!) then
$$(\Delta h_N)^n\to (\Delta h)^n$$
in the sense of measures. \footnote{A bit more precisely, we can
trivialize the line bundle $\det\ce_0^*[-2n]_\RR$ somehow and
consider $(\Delta h_N)^n$ and $(\Delta h)^n$ as measures on $M$. For
the very precise statement we refer to Section \ref{S:psh}.}

The proof of this result uses, among other things, the following new
multiplicativity property of the Baston operators:
\begin{eqnarray}\label{E:mult-bas}
\Delta(\ome\wedge\Delta \eta)=\Delta\ome\wedge\Delta \eta.
\end{eqnarray}
Here the operator $\Delta$ is understood on more general bundles
$\wedge^{k}\ce^*_0[-k]_\RR$ for various $k$'s, and $\ome,\eta$ are
sections of such more general bundles. This identity should be
compared with its obvious analogue from the case of complex
manifolds: for any differential forms $\ome,\eta$ one has
$$dd^c(\ome\wedge dd^c\eta)=dd^c\ome\wedge dd^c\eta.$$
The equality (\ref{E:mult-bas}) is proved in Section
\ref{S:baston-mult}. The proof uses general multiplicative structure
in spectral sequences; this is due to the fact that $\Delta$ itself
is defined with the use of spectral sequences.

\begin{remark}
The operators $\Delta$ on various bundles over $M$ mentioned
above fit into complexes of differential operators whose cohomology
can be interpreted as cohomology of the twistor space of $M$. On the
flat manifolds of any dimension some of these complexes with such an
interpretation were first constructed by Gindikin and Henkin
\cite{gindikin-henkin}, and the remaining ones by Henkin and
Polyakov \cite{henkin-polyakov}.
\end{remark}

\hfill

{\bf Acknowledgements.} I am grateful to G. Henkin for very useful
discussions; in particular some years ago he has explained to me the
relation between the Penrose transform and pluriharmonic functions
on the flat quaternionic space. I thank also A. \v Cap, V.
Palamodov, A. Swann, M. Verbitsky for very useful discussions
regarding quaternionic manifolds. I thank A. Beilinson, J.
Bernstein, V. Hinich, and Y. Varshavsky for very useful discussions
on the homological algebra. I thank also M. Verbitsky for numerous
remarks on the paper.

\section{Quaternionic manifolds.}\label{Ss:quatenionic_mflds}
In this section we remind the definition and basic properties of a
quaternionic manifolds due to S. Salamon \cite{salamon-82}. First
let us fix notation. Let $GL_n(\HH)$ denote the group of invertible
$n\times n$-matrices with quaternionic entries.  The group
$GL_n(\HH)\times GL_1(\HH)$ acts on $\HH^n$ by
\begin{eqnarray}\label{qm1}
(A,q)(x)=Axq^{-1}
\end{eqnarray}
where $(A,q)\in GL_n(\HH)\times GL_1(\HH),$ and $x\in \HH^n$ is a
column of quaternions. Thus we get a group homomorphism
\begin{eqnarray}\label{E:map-groups}
GL_n(\HH)\times GL_1(\HH)\to GL_{4n}(\RR).
\end{eqnarray}
The image of it is a closed subgroup denoted by $GL_n(\HH)\cdot
GL_1(\HH)$. Clearly the kernel is the multiplicative group of real
numbers $\RR^*$ imbedded diagonally into $GL_n(\HH)\times
GL_1(\HH)$. Let us denote by $\cg$ the subgroup of $GL_n(\HH)\times
GL_1(\HH)$ defined by
\begin{eqnarray}\label{D:G-grp}
\cg:=\{(A,q)\in GL_n(\HH)\times GL_1(\HH)|\, \det A\cdot |q|=1\}
\end{eqnarray}
where $\det$ of a quaternionic matrix is taken in the sense of
Diedonn\'e \footnote{See e.g. \cite{aslaksen},
\cite{alesker-bsm-03}; let us only mention that it is the only
homomorphism of groups $GL_n(\HH)\to \RR_{>0}^*$ whose restriction
to the subgroup of complex $n\times n$ matrices is equal to the
absolute value of the usual determinant of complex matrices.}, and
$|q|$ denotes the absolute value of the quaternion $q$.
Alternatively $\cg$ can be characterized as a connected subgroup
whose Lie algebra is equal to $\{(X,Y)\in gl_n(\HH)\times
gl_1(\HH)|\, \sum_{i=1}^nRe(X_{ii})+ Re(Y)=0\}$. The homomorphism
$$\cg\to GL_n(\HH)\cdot GL_1(\HH)$$
given by the composition of the imbedding $\cg\inj GL_n(\HH)\times
GL_1(\HH)$ with the homomorphism (\ref{E:map-groups}) is onto, and
the kernel is equal to $\{\pm 1\}$ imbedded diagonally.

\hfill

Let $M$ be a real smooth manifold of dimension $4n$. A {\itshape
quaternionic structure} on $M$ is a reduction of the structure group
of the tangent bundle of $M$ to $GL_n(\HH)\cdot GL_1(\HH)$ such that
there exists a torsion free $\gh$-connection on the tangent bundle
$TM$. The manifold $M$ with a quaternionic structure is called a
{\itshape quaternionic manifold}. Through the rest of the paper we
assume that $n>1$ unless otherwise stated; the case $n=1$ is rather
exceptional and requires a separate treatment.

\begin{remark}\label{qm2}
On a quaternionic manifold a torsion free connection as above is
{\itshape not} unique.
\end{remark}

Let us make an additional assumption about the quaternionic manifold
$M$ that the structure group of $TM$ has been lifted from $\gh$ to
$\cg$ (this is possible provided e.g. $H^2(M,\ZZ/2\ZZ)=0$; see
\cite{salamon-2Glob_Riem_Geom}). Let us fix such a lifting. This
assumption and this choice will influence none of the main results
of the paper since locally this always can be done, but then it can
be shown that all the main results are independent of these local
choices and make sense globally. Let $G_0\to M$ be the corresponding
principal $\cg$-bundle. Then the tangent bundle $TM$ is isomorphic
to $G_0\times _{\cg}\HH^n$ where the group $\cg$ acts on $\HH^n$ by
(\ref{qm1}). Let us define two other bundles over $M$
\begin{eqnarray}\label{qm3}
\ce_0=G_0\times_{\cg}\HH^n
\end{eqnarray}
where now the action of $\cg$ on $\HH^n$ is given by $(A,q)(x)=Ax$;
and the bundle
\begin{eqnarray}\label{qm4}
\ch_0=G_0\times_{\cg} \HH
\end{eqnarray}
where the action of $\cg$ on $\HH$ is given by $(A,q)(x)=xq^{-1}$.

It is easy to see that $\ce_0$ is a bundle of {\itshape right}
$\HH$-vector spaces of rank $n$ over $\HH$, while $\ch_0$ is a
bundle of {\itshape left} $\HH$-vector spaces of rank 1 over $\HH$.
Moreover we have a natural isomorphism
\begin{eqnarray}\label{qm5}
TM\simeq \ce_0\otimes_\HH \ch_0
\end{eqnarray}
(notice that the tensor product is over $\HH$).

\hfill

Let us study now the complexified tangent bundle $\cct
M:=TM\otimes_\RR \CC$. By (\ref{qm5})
\begin{eqnarray}\label{qm6}
\cct M=(\ce_0\otimes_\RR\CC)\otimes_{\cchh}(\ch_0\otimes_\RR\CC)
\end{eqnarray}
where $\cchh:=\HH\otimes_\RR\CC$ is a $\CC$-algebra
(non-canonically) isomorphic to the algebra $M_2(\CC)$ of $2\times
2$ complex matrices.

\begin{lemma}\label{qmL:1}
There exists an isomorphism of vector bundles
$$\cct M\simeq \ce_0\otimes_\CC\ch_0$$
where the tensor product is over $\CC$, and $\ce_0$, $\ch_0$ are
considered as $\CC$-vector bundles via the imbedding of fields
$\CC\inj \HH$ given by $\sqrt{-1}\mapsto I$.
\end{lemma}
\begin{remark}
The isomorphism constructed in the proof of the lemma will be used
in the rest of the paper. This isomorphism is almost canonical: it
depends on some non-canonical universal choice on the level of
linear algebra.
\end{remark}

{\bf Proof of Lemma \ref{qmL:1}.} Recall that $\cchh\simeq M_2(\CC)$
is a central simple $\CC$-algebra. It has a unique up to
(non-canonical) isomorphism simple right $\cchh$-module $T$. We
choose $T$ to be the right $\cchh$-submodule of $\cchh$ as follows
$$T:=\{x-\sqrt{-1}I\cdot x|\, x\in \HH\}=\{z\in \cchh|\, I\cdot
z=\sqrt{-1}z\}.$$ Let $E$ be a right $\HH$-vector space. We have a
functorial isomorphism of right $\cchh$-modules
\begin{eqnarray}\label{qm7}
E\otimes_\CC T\tilde\to E\otimes_\RR\CC
\end{eqnarray}
given by
\begin{eqnarray}\label{qm8}
\xi\otimes (x-\sqrt{-1}I x)\mapsto \xi\cdot x-\sqrt{-1}\xi\cdot
(Ix).
\end{eqnarray}
It is an easy exercise to check that this defines a well defined
morphism of right $\cchh$-modules $E\otimes_\CC T\to E\otimes_\RR
\CC$, and it is an isomorphism.

Similarly let $T'$ be the simple left $\cchh$-submodule of $\cchh$
$$T'=\{x-\sqrt{-1}xI|\, x\in \HH\}=\{z\in \cchh|\,
z\cdot I=\sqrt{-1}z\}.$$

For any left $\HH$-module $H$ we have the functorial isomorphism of
left $\cchh$-modules
\begin{eqnarray}\label{qm9}
T'\otimes_\CC H\tilde\to H\otimes_\RR\CC
\end{eqnarray}
given by
\begin{eqnarray}\label{qm10}
(x-\sqrt{-1}xI)\otimes h\mapsto xh-\sqrt{-1}xIh.
\end{eqnarray}
Then by (\ref{qm6}), (\ref{qm7}), (\ref{qm9}) we have
$$\cct M\simeq (\ce_0\otimes_\CC T)\otimes_{\cchh}(T'\otimes_\CC
\ch_0)\simeq(\ce_0\otimes_\CC\ch_0)\otimes_\CC(T\otimes_{\cchh}T').$$
it remains to observe that $T\otimes_{\cchh} T'\simeq \CC$. \qed

\hfill

It is well known that for $n>1$ the manifold $M^{4n}$ and the
bundles $\ce_0, \ch_0$ are {\itshape real analytic}. We denote by
$X$ a small enough complexification of $M$. This $X$ is a complex
manifold of complex dimension $4n$. We denote by $\ce, \ch$ the
holomorphic vector bundles over $X$ extending the above
$\ce_0,\ch_0$, i.e. $\ce|_M=\ce_0,\, \ch|_M=\ch_0$. Then there is an
isomorphism of holomorphic vector bundles $TX\simeq \ce\otimes
_\CC\ch$ extending the isomorphism $\cct M\simeq \ce_0\otimes_\CC
\ch_0$. Clearly $X$ is equipped with the complex conjugation
diffeomorphism
$$\sigma\colon X\to X$$
which is an anti-holomorphic involution, and the set of fixed points
is $X^\sigma =M$.

\hfill

Recall that $\ce_0$ is a bundle of right $\HH$-vector spaces, in
particular the multiplication by $J\in \HH$ on the right is an
$I$-anti-linear operator on fibers of $\ce_0$. By analytic
continuation (since everything is real analytic) we get the
following structure on $\ce$: on the total space of $\ce$, which is
a complex manifold, we have an anti-holomorphic map
$$\hat J\colon \ce\to \ce$$ such that for any $x\in M\subset X$ the
restriction of $\hat J$ to the fiber $\ce|_x$ coincides with the
right multiplication by $J$, and for any $z\in X$
\begin{eqnarray*}
\hat J\colon \ce|_z\tilde\to \ce|_{\sigma(z)}\mbox{ is
anti-linear}\label{qm11},\\
\hat J^2=-Id_{\ce|_z}\label{qm12}.
\end{eqnarray*}
Moreover $\hat J$ preserves the class of holomorphic sections of
$\ce$. This action on holomorphic sections of $\ce$ will be also
denoted by $\hat J$. Now we remind a definition (see e.g.
\cite{manin-gauge}, Ch. 2, \S 1.7).
\begin{definition}\label{D:real-quatern-str}
Let $X$ be a complex manifold.

(1) A {\itshape real structure} on $X$ is an antiholomorphic
involution $\sigma\colon X\to X$ (i.e. $\sigma^2=Id_X$). Clearly
$\sigma$ defines the anti-linear map on sections of the sheaf
$\co_X$ of holomorphic functions given by
$$f\mapsto \sigma(f):=\overline{f\circ \sigma}.$$

(2) Let $\cf$ be a coherent sheaf on $X$. Let $\sigma$ be a real
structure on $X$. A {\itshape real structure} on $\cf$ is an
antilinear map $\rho$ on sections of $\cf$ extending $\sigma$ from
$X$ and $\co_X$ (in particular $\rho(f\xi)=\sigma(f)\rho(\xi)$ for
any section $f$ of $\co_X$ and any section $\xi$ of $\cf$) such that
$\rho^2=Id$.

A {\itshape quaternionic structure} on $\cf$ is an antilinear map
$\rho$ on sections of $\cf$ extending $\sigma$ from $X$ and $\co_X$
such that $\rho^2=-Id$.
\end{definition}

The previous discussion implies that in our situation of a
quaternionic manifold the sheaf of holomorphic sections of $\ce$ is
equipped with the quaternionic structure $\hat J$.

Similarly the left multiplication by $J$ on $\ch_0$ induces an
anti-holomorphic diffeomorphism $\hat J\colon \ch\to \ch$ with the
similar properties. It follows that $\sigma,\hat J$ induce a
quaternionic structure on $\ch$. Taking the dual bundle we get a
quaternionic structure on $\ch^*$ which we will denote again by
$\hat J$. Notice also that it induces a real structure on the
projectivization $\PP(\ch^*)$, and the natural projection
$\PP(\ch^*)\to X$ commutes with the real structures on the two
spaces.

\hfill

Let us return back again to a general complex manifold $X$ with a
real structure $\sigma$. If $\cf_1$ and $\cf_2$ are coherent sheaves
on $X$ with either real or quaternionic structures $\rho_1$ and
$\rho_2$ respectively extending a real structure $\sigma$ on $X$,
then $\cf_1\otimes_{\co_X}\cf_2$ is equipped with $\rho_1\otimes
\rho_2$ (also extending $\sigma$) which is real if $\rho_1,\rho_2$
have the same type (i.e. real or quaternionic simultaneously), and
quaternionic if $\rho_1,\rho_2$ have different type. In particular
it follows that if $\cf$ is the sheaf of holomorphic sections of a
holomorphic vector bundle, and $\cf$ is equipped with a quaternionic
structure $\rho$, then $\wedge^k\cf$ is equipped with $\wedge^k\rho$
which is real for even $k$, and quaternionic for odd $k$.

Next, let $\cf$ be the sheaf of holomorphic sections of a
holomorphic vector bundle $\cv$ over $X$ with a real structure
$\rho$ extending $\sigma$. Let $X^\sigma$ be the set of fixed points
of $\sigma$. Then $X^\sigma$ is a real analytic submanifold of $X$.
Moreover if it is non-empty then $\dim_\RR X^\sigma=\dim_\CC X$, and
$X$ is a complexification of $X^\sigma$. Let us denote by $\cv_0$
the restriction of the vector bundle $\cv$ to $X^\sigma$. This
complex vector bundle $\cv_0$ is equipped with the fiberwise real
structure. The latter condition means that
$\cv_0=(\cv_0)_\RR\otimes_\RR \CC$ where $(\cv_0)_\RR$ is the real
vector bundle whose fiber over any point from $X^\sigma$ consists of
$\rho$-invariant vectors in the fiber of $\cv_0$ over that point.

\hfill

The above discussion implies that for a quaternionic manifold $M$
with a complexification $X$ and vector bundles $\ce,\ch$ as
previously, we obtain that each vector bundle
$\wedge^{2k}\ce^*\otimes (\det \ch^*)^{\otimes 2k}$, which also will
be denoted briefly as $\wedge^{2k}\ce^*[-2k]$, is equipped with a
real structure. Consequently its restriction
$\wedge^{2k}\ce_0^*\otimes (\det \ch_0^*)^{\otimes
2k}=:\wedge^{2k}\ce_0^*[-2k]$ to $M$ has a pointwise real structure,
namely
$$\wedge^{2k}\ce_0^*[-2k]=\wedge^{2k}\ce_0^*[-2k]_\RR\otimes_\RR\CC,$$
where $\wedge^{2k}\ce_0^*[-2k]_\RR$ is a real analytic vector bundle
over $M$ defined by the previous general construction.

These real vector bundles $\wedge^{2k}\ce_0^*[-2k]_\RR$ will be used
below, in particular in Section \ref{S:positive-currents}, where we
will introduce the notion of positivity of sections of these bundles
necessary to develop quaternionic pluripotential theory.

\section{The twistor space.}\label{Ss:twistor} For a quaternionic
manifold $M^{4n},\, n>1$, S. Salamon \cite{salamon-82} has
constructed a complex manifold $Z$ of complex dimension $2n+1$
called the {\itshape twistor space} of $M$. We remind now this
construction since it will be important in Section
\ref{S:penrose-baston} for the description of Baston's differential
operators.

We will assume that the vector bundles $\ch_0,\ce_0$ from Section
\ref{Ss:quatenionic_mflds} are defined globally over $M$ though this
assumption is not necessary; the whole construction can be done
first locally on $M$, and then it can be easily shown that it is
independent of local choices.

As a smooth manifold the twistor space
$$Z:=\PP(\ch_0^*)$$
is the (complex) projectivization of $\ch_0^*$. To describe the
complex structure on $Z$ it will be convenient to describe first the
complex structure on the total space of $\ch_0^*$ with the zero
section removed (in fact this space carries not only complex
structure, but a hypercomplex structure, see
\cite{pedersen-poon-swann}).


Let us choose a torsion free $\cg$-connection $\nabla$. For any
point $z\in \ch_0^*\backslash\{0\}$ the connection $\nabla$ induces
a decomposition of the tangent space $T_z\ch_0^*$ to the direct sum
of the vertical $\cv_z$ and the horizonal $\cl_z$ subspaces with
respect to the natural projection $p\colon \ch_0^*\to M$:
$$T_z\ch_0^*=\cv_z\oplus \cl_z.$$
Clearly $\cv_z=\ch_0^*|_{p(z)}$. Next the differential $dp\colon
\cl_z\to T_{p(z)}M$ is an isomorphism of $\RR$-vector spaces. One
equips $\cv_z=\ch_0^*|_{p(z)}$ with the complex structure $I$
(recall that $\ch_0^*$ is a quaternionic vector bundle). Let us
describe the complex structure on $\cl_z\tilde\to T_{p(z)}M$. We
follow \cite{salamon-2Glob_Riem_Geom}. Let $l\subset
\ch_0^*|_{p(z)}$ be the complex line spanned by $z$. The complex
structure on $T_{p(z)}M$ is uniquely characterized by the property
that the space of $(1,0)$-forms at $p(z)$ is equal to the subspace
$$\ce^*_0|_{p(z)}\otimes_\CC l\subset \ce_0^*|_{p(z)}\otimes_\CC
\ch_0^*|_{p(z)}\overset{\mbox{Lemma } \ref{qmL:1}}\simeq
T^*_{p(z)}M\otimes_\RR\CC.$$

Thus we got an almost complex structure on $\ch_0^*$. It is was
shown in \cite{pedersen-poon-swann}, Theorem 3.2, that it is
integrable and is independent of a choice of a torsion free
connection $\nabla$.

Now the non-zero complex numbers $\CC^*$ act holomorphically by the
product on $\ch_0^*\backslash\{0\}$, and the quotient is equal to
the twistor space $Z=\PP(\ch_0^*)$. Hence $Z$ carries a complex
structure.

Observe moreover that we have constructed also a holomorphic
principal $\CC^*$-bundle $\ch_0^*\backslash\{0\}\to Z$. Let
$\co_Z(-1)$ be the holomorphic line bundle over $Z$ defined by
$$\co_Z(-1):=\ch_0^*\backslash\{0\}\times_{\CC^*}\CC$$
where $\CC^*$ acts on $\CC$ by the usual multiplication. The dual
line bundle is denoted by $\co_Z(1)$; this holomorphic line bundle
is called sometimes the Swann bundle. Notice also that as a
{\itshape smooth} complex line bundle $\co_Z(1)$ is easily described
as follows: it is usual dual Hopf line bundle over $Z=\PP(\ch^*_0)$.
Let us emphasize that $\co_Z(1)$ is defined globally only under the
assumption that the structure group of $TM$ is lifted to $\cg$.
However $Z$ itself and $\co_Z(2):=(\co_Z(1))^{\otimes 2}$ are
defined globally independently of any local choices. In this paper
we will really use only $\co_Z(2)$ and its tensor powers.


\hfill

As in Section \ref{Ss:quatenionic_mflds} we denote by $X$ a small
enough complexification of $M$, and by $\ch$ the vector bundle over
$X$ which is the (unique) holomorphic extension of the real analytic
vector bundle $\ch_0$ from $M$ to $X$. Let us consider the complex
analytic manifold $\PP(\ch^*)$. We have the obvious holomorphic map
$\tau\colon \PP(\ch^*)\to X$. It is obvious that
$$\tau^{-1}(M)=\PP(\ch_0^*)=Z.$$
(Warning: $Z=\tau^{-1}(M)$ is {\itshape not} a complex submanifold
of $\PP(\ch^*)$.)

The next non-trivial claim is that there exists a {\itshape
holomorphic} submersion
$$\eta\colon \PP(\ch^*)\to Z$$
which is uniquely characterized by the property that
$\eta|_{\tau^{-1}(M)=Z}=Id_Z$ (see \cite{baston-92},
\cite{cap-slovak-soucek-II}).

Thus we get a diagram of holomorphic submersions
\begin{eqnarray}\label{Diag:twistor}
Z\overset{\eta}{\leftarrow}\PP(\ch^*)\overset{\tau}{\to}X.
\end{eqnarray}

Let us assume again that the structure group of $TM$ is lifted to
$\cg$. Over the complex manifold $\PP(\ch^*)$ we have the usual
holomorphic dual Hopf line bundle which will be denoted by
$\tilde\co(1)$. Observe that its restriction to $Z=\tau^{-1}(M)$ is
equal to $\co_Z(1)$.

\begin{lemma}\label{L:iso-bundles-ex-un}
There exists a unique isomorphism of holomorphic line bundles
\begin{eqnarray}\label{Iso-oo}
\Phi\colon \eta^*(\co_Z(1))\tilde\to \tilde\co(1)
\end{eqnarray}
such that the restriction of $\Phi$ to $Z=\tau^{-1}(M)$ is equal to
the identity.
\end{lemma}
{\bf Proof.} First let us prove the uniqueness; we will show that in
fact $\Phi$ is unique locally on $\PP(\ch^*)$, more precisely $\Phi$
is unique in a neighborhood of any point $p\in \tau^{-1}(M)$. Any
such point $p$ has a neighborhood and holomorphic coordinates
$(z_1,\dots,z_{4n},w)$ such that
$$\tau^{-1}(M)=\{(z_1,\dots,z_{4n},w)|\, Im(z_1)=\dots
=Im(z_{4n})=0\}.$$ Assume that we have two such isomorphisms $\Phi$
and $\Phi'$. Then there exists a holomorphic function $f$ such that
$$\Phi'=f\cdot \Phi$$ and $f|_{\tau^{-1}(M)}=1$. In order to show
that $f=1$ identically, let us decompose $f$ into the holomorphic
Taylor series:
$$f=\sum_{\alp,b} f_\alp z^\alp w^b$$
where $\alp=(\alp_1,\dots,\alp_{4n})\in \ZZ_{\geq 0}^{4n},\, b\in
\ZZ_{\geq 0},\, z=(z_1,\dots, z_{4n})$. We know that $f=1$ whenever
$Im(z_1)=\dots =Im(z_{4n})=0$. It follows that all partial
derivatives of $f$ of positive degree vanish. Hence $f=1$
identically, and $\Phi'=\Phi$.

\hfill

Now let us prove the existence of $\Phi$. Due to the uniqueness of
such an isomorphism, which was proved even locally on $\PP(\ch^*)$,
it suffices to prove the existence locally in a neighborhood of an
arbitrary point $p\in \tau^{-1}(M)$; here one should make $X$
smaller if necessary and use the properness of the map $\tau$.

The obvious (identity) isomorphism between the vector bundles
$\eta^*\co_Z(1)|_{\tau^{-1}(M)}$ and $\tilde \co(1)|_{\tau^{-1}(M)}$
over $\tau^{-1}(M)$ is real analytic (since all the manifolds,
morphisms and vector bundles are real analytic). Let us fix
holomorphic trivializations of $\eta^*\co_Z(1)$ and $\tilde \co(1)$
in a neighborhood $U\subset X$ of $p$. Then our real analytic
isomorphism over $\tau^{-1}(M)$ between the trivialized line bundles
is given by a non-vanishing real analytic function
$$g\colon \tau^{-1}(M)\cap U\to \CC.$$

Let us decompose $g$ into (real) Taylor series converging in
$\tau^{-1}(M)\cap U$
\begin{eqnarray}\label{E:COR-series}
g=\sum_{\alp,b,c} g_{\alp,b,c} x^\alp w^b \bar w^c
\end{eqnarray}
where $x=(x_1,\dots,x_{4n})\in \RR^{4n}$ with $x_i=Re (z_i)$,
$\alp\in \ZZ_{\geq 0}^{4n},b,c\in \ZZ_{\geq 0}$.

Next it is clear that the restriction of $g$ to any complex curve
$\tau^{-1}(m)\cap U$, with $m\in M$ being an arbitrary point, is a
holomorphic function (this is because
$\eta^*\co_Z(1)|_{\tau^{-1}(m)}\simeq \tilde\co(1)|_{\tau^{-1}(m)}$
is the Hopf bundle over $\tau^{-1}(m)\simeq \CC\PP^1$).

This is equivalent to say that for any fixed $x\in \RR^{4n}$, $g$ is
a holomorphic function in $w$. This implies that actually $\bar w$
does not appear in (\ref{E:COR-series}):
$$g=\sum_{\alp,b} g_{\alp,b} x^\alp w^b.$$
Since the Taylor series converges in $U$, it is a restriction to
$\tau^{-1}(M)$ of the holomorphic function
$$\tilde g:=\sum_{\alp,b} g_{\alp,b} z^\alp w^b.$$
Then $\tilde g$ induces a local isomorphism
$$\eta^*\co_Z(1)|_U\tilde\to \tilde\co(1)|_U$$
which is equal to the identity on $\tau^{-1}(M)\cap U$. \qed

\hfill

Let us discuss now quaternionic structures on the above line
bundles. In Section \ref{Ss:quatenionic_mflds} we have defined a
real structure on $X$ which is just the complex conjugation
$\sigma\colon X\to X$. Also we have defined a quaternionic structure
$\hat J\colon \ch^*\to \ch^*$. It induces a real structure
$\tilde\sigma$ on the projectivization of $\ch^*$:
$$\tilde\sigma\colon \PP(\ch^*)\to \PP(\ch^*).$$
Moreover $\hat J$ induces a quaternionic structure $\tilde J$ on the
line bundle $\tilde \co(1)$ which extends $\tilde\sigma$ (see
Definition \ref{D:real-quatern-str}):
$$\tilde J(\xi)=\xi\cdot \hat J$$
(notice that $\hat J$ acts on the right since $\ch^*_0$ is a bundle
of right $\HH$-modules).

The restriction of $\tilde \sigma$ to $Z=\tau^{-1}(M)\subset
\PP(\ch^*)$ is antiholomorphic, and hence it can be considered as a
real structure on $Z$ which will be denoted by $\sigma_Z$. The map
$\eta\colon \PP(\ch^*)\to Z$ intertwines $\tilde\sigma$ and
$\sigma_Z$:
\begin{eqnarray}\label{E:sigma-sigma}
\eta\circ \tilde\sigma=\sigma_Z\circ \eta.
\end{eqnarray}
Indeed $\sigma^{-1}_Z\circ \eta\circ\tilde\sigma\colon \PP(\ch^*)\to
Z$ is a holomorphic map whose restriction to $Z=\tau^{-1}(M)$ is
equal to $Id_Z$. Hence this map must be equal to $\eta$.

\hfill

Next, the restriction of $\tilde J$ to $\co_Z(1)$ induces a
quaternionic structure on $\co_Z(1)$ which we denote by  $J_Z$.

The holomorphic line bundle $\eta^*(\co_Z(1))$ has an induced
quaternionic structure $\eta^*J_Z$ which extends the real structure
$\tilde\sigma$ because of (\ref{E:sigma-sigma}). The isomorphism
$\Phi\colon \eta^*(\co_Z(1))\tilde\to \tilde\co(1)$ from
(\ref{Iso-oo}) is compatible with quaternionic structures:
$$\Phi\circ \eta^*J_Z=\tilde J\circ \Phi.$$
Indeed $\tilde J^{-1}\circ \Phi\circ \eta^* J_Z\colon
\eta^*(\co_Z(1))\tilde\to \tilde\co(1)$ is a holomorphic isomorphism
of holomorphic vector bundles whose restriction to $Z=\tau^{-1}(M)$
is equal to identity.

\section{The Penrose transform and the Baston
complexes.}\label{S:penrose-baston} The goal of this section is to
describe a construction due to Baston \cite{baston-92} of certain
complexes of vector bundles over a quarternionic manifold $M^{4n}$,
$n>1$, where the differentials are differential operators of either
first or second order. These complexes depend only on the
quaternionic structure of $M$, in particular they are equivariant
under quaternionic automorphisms of $M$. One of the differential
operators of the second order in one of these complexes will be
necessary in our definition of plurisubharmonic functions in Section
\ref{S:psh} below. Some of the other operators will be useful for
technical reasons, e.g. in the proof of a quaternionic
generalization of theorems of Aleksandrov and
Chern-Levine-Nirenberg.

The Baston approach \cite{baston-92} to construct these complexes is
based on the use of the Penrose transform. Thus we will have to
remind this notion in this section. It is convenient to assume
existence and fix a lifting of the structure group of $TM$ to
$\cg\subset \gth$ defined in Section \ref{Ss:quatenionic_mflds} (in
applications to the pluripotential theory a lifting always exists
locally). Global results are independent of such local liftings.

We keep the notation of Sections \ref{Ss:quatenionic_mflds},
\ref{Ss:twistor}. But for brevity we will denote in this section
$F:=\PP(\ch^*)$. Consider again as in Section \ref{Ss:twistor} the
holomorphic maps
$$Z\overset{\eta}{\leftarrow} F\overset{\tau}{\to} X$$
(recall that $X$ is a small enough complexification of $M$, and $Z$
is the twistor space). Let $\cl$ be a holomorphic vector bundle over
$Z$; by the abuse of notation we will denote also by $\cl$ the sheaf
of holomorphic sections of it. Let $\eta^{-1}\cl$ denote the
pull-back of $\cl$ under the map $\eta$ in the category of abstract
sheaves (thus in particular $\eta^{-1}\cl$ is {\itshape not} a sheaf
of holomorphic sections of anything). The {\itshape Penrose
transform} of $\cl$, by the definition, is $R\tau_*(\eta^{-1}\cl)$
where $R\tau_*$ is the push-forward morphism under $\tau$ between
the derived categories of sheaves:
$$R\tau_*\colon D^+(Sh_F)\to D^+(Sh_X)$$
where $D^+(Sh_X)$ denotes the bounded from below derived category of
sheaves of $\CC$-vector spaces on $X$.

\hfill

Let $(\Omega^\bullet_{F/Z},d)$ denote the complex of holomorphic
relative differential forms with respect to the map $\eta\colon F\to
Z$. We have the resolution $\Ome^\bullet(\cl)$ of the sheaf
$\eta^{-1}\cl$ by coherent sheaves
\begin{eqnarray}\label{E:L-resolution}
0\to\eta^{-1}\cl\to
\eta^*\cl\overset{d}{\to}\eta^*\cl\otimes_{\co_F}\Ome^1_{F/Z}
\overset{d}{\to}\eta^*\cl\otimes_{\co_F}\Ome^2_{F/Z}\overset{d}{\to}\dots
\end{eqnarray}
where $\eta^*\cl$ denotes the pull-back of $\cl$ in the category of
quasi-coherent sheaves, all the tensor products are over the sheaf
$\co_F$ of holomorphic functions on $F$, and the differentials are
the holomorphic de Rham differentials. Hence the Penrose transform
$R\tau_*(\eta^{-1}\cl)$ is equal to
\begin{eqnarray}\label{E:Penrose-complex}
R\tau_*(\Ome^\bullet(\cl))=R\tau_*\left(0\to
\eta^*\cl\overset{d}{\to}\eta^*\cl\otimes_{\co_F}\Ome^1_{F/Z}
\overset{d}{\to}\eta^*\cl\otimes_{\co_F}\Ome^2_{F/Z}\overset{d}{\to}\dots\right).
\end{eqnarray}

Let us consider the hypercohomology spectral sequence for
(\ref{E:Penrose-complex}) such that its first terms, denoted by
$E_1^{p,q}(\cl)$, are
$$E_1^{p,q}(\cl)=R^q\tau_*(\eta^*\cl\otimes_{\co_F}\Ome^p_{F/Z}).$$
Since the fibers of $\tau$ are complex projective lines and all the
sheaves $\eta^*\cl\otimes_{\co_F}\Ome^p_{F/Z}$ are coherent, we have
$$E_1^{p,q}(\cl)=0 \mbox{ for } q\ne 0,1.$$

Next assume that $\cl$ has real (resp. quaternionic) structure
$\rho$ extending the real structure $\sigma_Z$ on $Z$ (see
Definition \ref{D:real-quatern-str} in Section
\ref{Ss:quatenionic_mflds}). Then clearly $\eta^*\cl$ has real
(resp. quaternionic) structure $\eta^*\rho$ extending the real
structure $\tilde\sigma$ on $F=\PP(\ch^*)$. The sheaves
$\eta^*\cl\otimes_{\co_F}\Ome^p_{F/Z}$ are equipped with a real
(resp. quaternionic) structure $\rho^p$ as follows:
$$\rho^p(\xi\otimes \ome)=(\eta^*\rho)(\xi)\otimes
\overline{\tilde\sigma^*(\ome)}$$ where the bar denotes the complex
conjugation on differential forms. It is easy to see that the
holomorphic relative de Rham differential (\ref{E:L-resolution}) is
compatible with $\rho^p$'s:
$$d\circ \rho^p=\rho^{p+1}\circ d.$$
Since $\tau\colon F\to X$ intertwines the real structures
$\tilde\sigma$ and $\sigma$, it follows that all the terms
$E^{p,q}_r(\cl)$ of our spectral sequence are equipped with real
(resp. quaternionic) structure, and the differentials $d_r^{p,q}$
commute with it.

\hfill

In order to construct the Baston complexes, Baston has chosen
$$\cl=\co_Z(-k):=\co_Z(-1)^{\otimes k},\, 2\leq k\leq 2n,$$
computed $E_1^{p,q}(\co_Z(-k))$, and appropriate differentials
$d_1^{p,q},d_2^{p,q}$. We will use these computations, so let us
describe them. For a sheaf $\cs$ on $X$ and an integer $l$ let us
denote for brevity
$$\cs[l]:=\cs\otimes_{\co_X}(\det\ch)^{\otimes l}.$$

\begin{proposition}[Baston \cite{baston-92}]\label{P:terms-spec-seq}
Let $k\geq 2$. Then one has
\begin{eqnarray*}
E_1^{p,q}(\co_Z(-k))=0 \mbox{ for } q\ne 0,1;\\
E_1^{p,0}(\co_Z(-k))=
\left\{\begin{array}{ccc}
        Sym^{p-k}\ch\otimes \wedge^p\ce^*[-p]
        &\mbox{ if } &p-k\geq 0,\\
        0& \mbox{ if } &p-k <0\end{array}\right. ;\\
E_1^{p,1}(\co_Z(-k))=\left\{\begin{array}{ccc}
                            Sym^{k-p-2}\ch^*\otimes \wedge^p\ce^*[-p-1]&\mbox{ if }&k-p-2\geq 0,\\
                            0&\mbox{ if }&k-p-2<0
                            \end{array}\right.
\end{eqnarray*}
where all the tensor products are over $\co_X$.
\end{proposition}
Notice that in particular
$$E_1^{k-1,q}(\co_Z(-k))=0 \mbox{ for any } q.$$
The proposition implies that the first terms of the spectral
sequence look as follows:
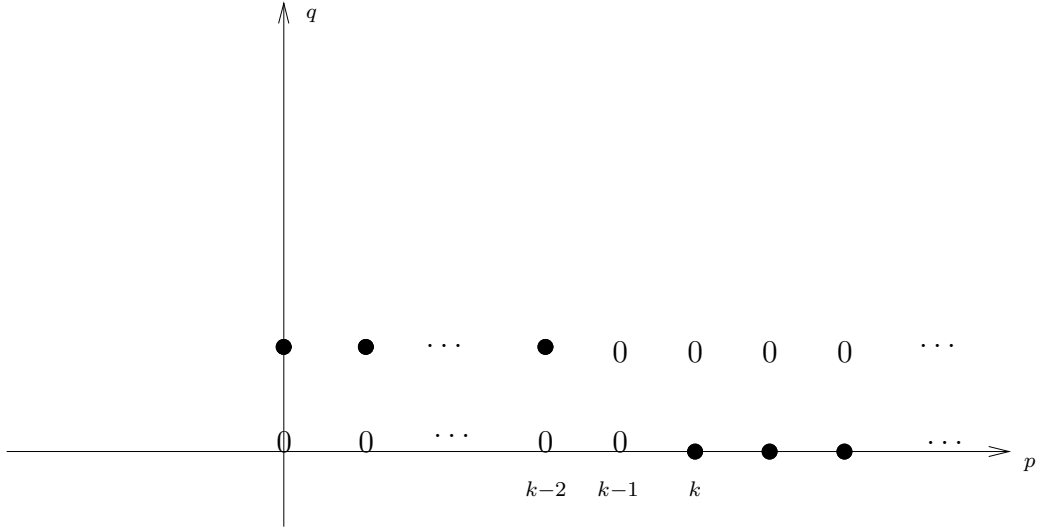
\begin{figure}[h]
\setlength{\unitlength}{0.00087489in}
\begingroup\makeatletter\ifx\SetFigFontNFSS\undefined%
\gdef\SetFigFontNFSS#1#2#3#4#5{%
  \reset@font\fontsize{#1}{#2pt}%
  \fontfamily{#3}\fontseries{#4}\fontshape{#5}%
  \selectfont}%
\fi\endgroup%
{\renewcommand{\dashlinestretch}{30}
\begin{picture}(6254,3210)(0,-10)
\put(1677,1092){\blacken\ellipse{90}{90}}
\put(1677,1092){\ellipse{90}{90}}
\put(2172,1092){\blacken\ellipse{90}{90}}
\put(2172,1092){\ellipse{90}{90}}
\put(3252,1092){\blacken\ellipse{90}{90}}
\put(3252,1092){\ellipse{90}{90}}
\put(4152,462){\blacken\ellipse{90}{90}}
\put(4152,462){\ellipse{90}{90}}
\put(4602,462){\blacken\ellipse{90}{90}}
\put(4602,462){\ellipse{90}{90}}
\put(5052,462){\blacken\ellipse{90}{90}}
\put(5052,462){\ellipse{90}{90}}
\path(1677,12)(1677,3162)
\path(1707.000,3042.000)(1677.000,3162.000)(1647.000,3042.000)
\path(12,462)(6042,462)
\path(5922.000,432.000)(6042.000,462.000)(5922.000,492.000)
\put(1632,462){\makebox(0,0)[lb]{\smash{{\SetFigFontNFSS{12}{14.4}{\rmdefault}{\mddefault}{\updefault}0}}}}
\put(2127,462){\makebox(0,0)[lb]{\smash{{\SetFigFontNFSS{12}{14.4}{\rmdefault}{\mddefault}{\updefault}0}}}}
\put(3207,462){\makebox(0,0)[lb]{\smash{{\SetFigFontNFSS{12}{14.4}{\rmdefault}{\mddefault}{\updefault}0}}}}
\put(2577,552){\makebox(0,0)[lb]{\smash{{\SetFigFontNFSS{12}{14.4}{\rmdefault}{\mddefault}{\updefault}$\ldots$}}}}
\put(3657,462){\makebox(0,0)[lb]{\smash{{\SetFigFontNFSS{12}{14.4}{\rmdefault}{\mddefault}{\updefault}0}}}}
\put(3657,1002){\makebox(0,0)[lb]{\smash{{\SetFigFontNFSS{12}{14.4}{\rmdefault}{\mddefault}{\updefault}0}}}}
\put(4557,1002){\makebox(0,0)[lb]{\smash{{\SetFigFontNFSS{12}{14.4}{\rmdefault}{\mddefault}{\updefault}0}}}}
\put(5007,1002){\makebox(0,0)[lb]{\smash{{\SetFigFontNFSS{12}{14.4}{\rmdefault}{\mddefault}{\updefault}0}}}}
\put(4107,1002){\makebox(0,0)[lb]{\smash{{\SetFigFontNFSS{12}{14.4}{\rmdefault}{\mddefault}{\updefault}0}}}}
\put(2532,1092){\makebox(0,0)[lb]{\smash{{\SetFigFontNFSS{12}{14.4}{\rmdefault}{\mddefault}{\updefault}$\ldots$}}}}
\put(5547,507){\makebox(0,0)[lb]{\smash{{\SetFigFontNFSS{12}{14.4}{\rmdefault}{\mddefault}{\updefault}$\ldots$}}}}
\put(5502,1092){\makebox(0,0)[lb]{\smash{{\SetFigFontNFSS{12}{14.4}{\rmdefault}{\mddefault}{\updefault}$\ldots$}}}}
\put(1812,3072){\makebox(0,0)[lb]{\smash{{\SetFigFontNFSS{12}{14.4}{\rmdefault}{\mddefault}{\updefault}$\scriptstyle q$}}}}
\put(6132,372){\makebox(0,0)[lb]{\smash{{\SetFigFontNFSS{12}{14.4}{\rmdefault}{\mddefault}{\updefault}$\scriptstyle p$}}}}
\put(3132,192){\makebox(0,0)[lb]{\smash{{\SetFigFontNFSS{12}{14.4}{\rmdefault}{\mddefault}{\updefault}$\scriptstyle k-2$}}}}
\put(3567,192){\makebox(0,0)[lb]{\smash{{\SetFigFontNFSS{12}{14.4}{\rmdefault}{\mddefault}{\updefault}$\scriptstyle k-1$}}}}
\put(4052,192){\makebox(0,0)[lb]{\smash{{\SetFigFontNFSS{12}{14.4}{\rmdefault}{\mddefault}{\updefault}
$\scriptstyle k$}}}}
\end{picture}
}
\caption{The first term $E^{p,q}_1({\mathcal O}_z(-k))$ of
the spectral sequence.}
\end{figure}

\newpage

The first differentials in this spectral sequence
$$d_1^{p,q}\colon E_1^{p,q}\to E_1^{p+1,q}$$
are differential operators of the first order (see
\cite{baston-92}).

The only non-zero differential in the {\itshape second} term
spectral sequence is
$$d_2^{k-2,1}\colon E_2^{k-2,1}\to E_2^{k,0}.$$
This clearly induces a morphism of sheaves
$$\Delta\colon \wedge^{k-2}\ce^*[-k+1]\to \wedge^k\ce^*[-k].$$
This $\Delta$ is a differential operator of the second order
\cite{baston-92}. Putting together all $d_1^{p,q}$ and $\Delta$ one
gets the Baston complex of coherent sheaves on $X$ for any $2\leq
k\leq 2n$:
\begin{eqnarray*}
0\to Sym^{k-2}\ch^*[-1]\overset{d_1^{0,1}}{\to}
Sym^{k-3}\ch^*\otimes\ce^*[-2]\overset{d_1^{1,1}}{\to}\dots\\
\dots\overset{d_1^{k-3,1}}{\to}\wedge^{k-2}\ce^*[-k+1]\overset{\Delta}{\to}\wedge^k\ce^*[-k]\overset{d_1^{k,0}}{\to}\\
\overset{d_1^{k,0}}{\to}\ch\otimes\wedge^{k+1}\ce^*[-k-1]\overset{d_1^{k+1,0}}{\to}Sym^2\ch\otimes
\wedge^{k+2}\ce^*[-k-2]\overset{d_1^{k+2,0}}{\to}\dots\\
\dots\overset{d_1^{2n-1,0}}{\to} Sym^{2n-k}\ch\otimes
\wedge^{2n}\ce^*[-2n]\to 0
\end{eqnarray*}
where all the tensor products are over $\co_X$. We call this complex
of coherent sheaves the {\itshape holomorphic Baston complex}.
Baston shows \cite{baston-92}, Proposition 10, that this complex of
holomorphic sheaves is acyclic except at the left. Moreover when the
quaternionic manifold $M$ is a small enough ball (not necessarily
flat as a quaternionic manifold), the complex of the global sections
of the Baston complex is a resolution of $H^1(Z,\co_Z(-k))$.

\hfill

Let us discuss now the quaternionic and the real structures on the
Baston complex. On each term of it let us take just the tensor
product of the quaternionic structures on $\ch,\ce$ and their duals
defined in Section \ref{Ss:quatenionic_mflds}. Then for each even
$k$ we get a real structure on each term of the complex (this case
will be particularly important for pluripotential theory), and for
odd $k$ we get a quaternionic structure. The differentials of the
Baston complex commute with these real (resp. quaternionic)
structures.

\hfill

Next passing to a germ of $M$ inside of $X$ we get the complex of
differential operators between the sheaves on $M$ of {\itshape real
analytic} sections of the corresponding vector bundles over $M$ (in
the above holomorphic Baston complex, one just replaces $\ch,\ce$
with $\ch_0,\ce_0$ respectively everywhere). The differential
operators in that complex have real analytic coefficients; we will
denote them by the same letters as in the holomorphic Baston
complex. They extend uniquely by continuity to morphisms of sheaves
of {\itshape infinitely smooth} sections of these vector bundles. In
particular we get the following complex:
\begin{eqnarray*}
0\to C^\infty(M,Sym^{k-2}\ch^*_0[-1])\overset{d_1^{0,1}}{\to}
C^\infty(M,Sym^{k-3}\ch^*_0\otimes\ce^*_0[-2])\overset{d_1^{1,1}}{\to}\dots\\
\dots\overset{d_1^{k-3,1}}{\to}C^\infty(M,\wedge^{k-2}\ce^*_0[-k+1])
\overset{\Delta}{\to}C^\infty(M,\wedge^k\ce^*_0[-k])\overset{d_1^{k,0}}{\to}\\
\overset{d_1^{k,0}}{\to}C^\infty(M,\ch_0\otimes\wedge^{k+1}\ce^*_0[-k-1])\overset{d_1^{k+1,0}}{\to}C^\infty(M,Sym^2\ch_0\otimes
\wedge^{k+2}\ce^*_0[-k-2])\overset{d_1^{k+2,0}}{\to}\dots\\
\dots\overset{d_1^{2n-1,0}}{\to}C^\infty(M, Sym^{2n-k}\ch_0\otimes
\wedge^{2n}\ce^*_0[-2n])\to 0.
\end{eqnarray*}
We will call this complex the {\itshape smooth Baston complex}.
\begin{remark}
Baston claims that the smooth Baston complex is also acyclic except
at the left. Moreover when $M$ is a small enough ball, the complex
of its global sections computes $H^1(Z,\co_Z(-k))$. We will not use
however these facts in this paper.
\end{remark}

\hfill

The case $k=2$ will be particularly important; let us write the
beginning of the smooth Baston complex in this case:
$$0\to C^\infty(M,\co_M[-1]):=C^\infty(M,\det\ch_0^*)\overset{\Delta}{\to}
C^\infty(M,\wedge^2\ce_0^*[-2])\overset{d_1^{2,0}}{\to}\dots.$$

\begin{remark}\label{R:normalization}
Later on in Section \ref{S:hypercomplex} we will change slightly the
definitions of $\Delta$ by multiplying it by an appropriate constant
in order to make it compatible with conventions in the flat case.
\end{remark}

As a side remark let us mention that some other complexes on
quaternionic manifolds were considered in \cite{widdows}.

\section{Multiplicative properties of the Baston
operators.}\label{S:baston-mult} Recall that on a quaternionic
manifold $M^{4n}$, $n>1$, we have the Baston differential operators
for any $2\leq k\leq 2n$
$$\Delta\colon C^{\infty}(M,\wedge^{k-2}\ce_0^*[-k+1])\to
C^\infty(M,\wedge^k\ce_0^*[-k]).$$ In this section we discuss the
relation of these operators to the obvious wedge product
$$\wedge^p\ce_0^*[i]\otimes\wedge^q\ce_0^*[j]\to
\wedge^{p+q}\ce_0^*[i+j].$$ This relation will be important in the
proof of the quaternionic generalization of the theorems of
Aleksandrov and Chern-Levine-Nirenberg. The main result of this
section is the following proposition which is a new result, to the
best of our knowledge.
\begin{proposition}\label{P:baston-mult}
Let $k,l\geq 2,\, k+l\leq 2n$. Then for any
$$\ome\in C^\infty(M,\wedge^{k-2}\ce_0^*[-k+1]),\, \eta\in
C^\infty(M,\wedge^{l-2}\ce_0^*[-l+1])$$ one has
\begin{eqnarray}\label{E:Delta-mult}
\Delta(\omega\wedge \Delta\eta)=\Delta\omega\wedge \Delta\eta.
\end{eqnarray}
\end{proposition}
{\bf Proof.} We may and will assume that $\ome$ and $\eta$ are real
analytic. Hence they can be considered as holomorphic sections over
$X$ of the bundles $\wedge^{k-2}\ce^*[-k+1]$ and
$\wedge^{l-2}\ce^*[-l+1]$ respectively. The result follows from
rather general properties of the multiplicative structure in
spectral sequences.

Recall that the (holomorphic) Baston operator was defined as the
only non-zero second differential in the hypercohomology spectral
sequence for the complex
$$\Ome^\bullet(\co_Z(-k))=\left(0\to\tilde\co(-k)\overset{d}{\to}\tilde\co(-k)\otimes_{\co_F}\Ome^1_{F/Z}
\overset{d}{\to}\tilde\co(-k)\otimes_{\co_F}\Ome^2_{F/Z}\overset{d}{\to}\dots
\right)$$ and similarly with $l$ instead of $k$.

First obviously $\co_Z(-k)\otimes_{\co_Z}\co_Z(-l)=\co_Z(-(k+l))$.
Hence we have the obvious natural morphism of complexes of sheaves
\begin{eqnarray}\label{E:mor-complexes}
\Ome^\bullet(\co_Z(-k))\otimes_{\underline{\CC}}\Ome^\bullet(\co(-l))\to
\Ome^\bullet(\co(-k-l))
\end{eqnarray}
where the tensor product is taken in the sense of complexes of
sheaves over the sheaf $\underline{\CC}$ of locally constant
$\CC$-valued functions. We will need few basic general facts on the
multiplicative structure in spectral sequences. We failed to find a
precise reference to these facts, but they seem to be a common
knowledge among experts in homological algebra. The author has
learned them from A. Beilinson \cite{beilinson}.

First let us equip the complexes $\Ome^\bullet(\co_Z(-k))$ and
$\Ome^\bullet(\co_Z(-l))$ with the stupid filtration, and to replace
them with a resolution in the bounded from below filtered derived
category $D^+F$ of sheaves of $\underline{\CC}$-modules (see e.g.
\cite{illusie}, Ch. V). The need to work with the filtered derived
category comes from the fact that in $D^+F$ we have both the derived
push-forward $R\tau_*$, tensor product, and the spectral sequence
related to $R\tau_*$. The usual derived category is not sufficient
since the spectral sequence is not defined on elements of it; on the
other hand in the category of actual complexes $R\tau_*$ is not well
defined.

There is a general notion of tensor product of abstract spectral
sequences \cite{mccleary}. Without giving all the formal details and
definitions, let us only explain the main idea. Given two spectral
sequences with terms and differentials $(E_r^{'p,q},d_r^{'p,q})$ and
$(E_r^{''p,q},d_r^{''p,q})$ respectively, the $E_r^{p,q}$-term of
their tensor product is defined by
$$E_r^{p,q}:=\bigoplus_{\begin{array}{c}
                         p'+p''=p\\
                         q'+q''=q
                         \end{array}}E_r^{'p',q'}\otimes_{\underline{\CC}}E_r^{''p'',q''}.$$
The differentials in the tensor product of the spectral sequences,
which will be denoted by $d^{p,q}_{\otimes, r}$, are defined as the
differential in the tensor product of complexes:
$$d_{\otimes,r}^{p,q}:=\bigoplus_{\begin{array}{c}
                         p'+p''=p\\
                         q'+q''=q
                         \end{array}}\left(d_r^{'p',q'}\otimes
                         1+(-1)^{p'+q'}\cdot 1\otimes
                         d_r^{''p'',q''}\right).$$
Then $(E_r^{p,q},d_{\otimes,r}^{p,q})$ is a new spectral sequence,
i.e. the cohomology of $d_{\otimes,r}$ is naturally isomorphic to
$E_{r+1}$. This is due to the fact that we work in the category of
sheaves over a field.

The key property of this tensor product of spectral sequences
computing $R\tau_*$ is that there exists a canonical morphism from
the tensor product of spectral sequences of two complexes of sheaves
of $\underline{\CC}$-modules to the spectral sequence (again
computing $R\tau_*$) of the tensor product over $\underline{\CC}$ of
these complexes. This fact and the morphism (\ref{E:mor-complexes})
lead to a canonical morphism of spectral sequences
\begin{eqnarray}\label{E:mor-spec-seq}
E(\co_Z(-k))\otimes E(\co_Z(-l))\to E(\co_Z(-k-l))
\end{eqnarray}
where the first expression denotes the tensor product of spectral
sequences mentioned above. In particular for $r=1$ we get morphism
$$E_1^{k-2,1}(\co_Z(-k))\otimes_{\underline{\CC}}E_1^{l,0}(\co_Z(-l))\to
E_1^{k+l-2,1}(\co_Z(-k-l)).$$ Computing these terms by Proposition
\ref{P:terms-spec-seq} we get a morphism of sheaves of
$\underline{\CC}$-modules
$$\wedge^{k-2}\ce^*[-k+1]\otimes_{\underline{\CC}}\wedge^l\ce^*[-l]\to
\wedge^{k+l-2}\ce^*[-k-l+1].$$ It is easy to see that this morphism
coincides with the composition of the natural morphism of sheaves
$$\wedge^{k-2}\ce^*[-k+1]\otimes_{\underline{\CC}}\wedge^l\ce^*[-l]\to
\wedge^{k-2}\ce^*[-k+1]\otimes_{\co_X}\wedge^l\ce^*[-l]$$ with the
wedge product.

By the definition of morphism of spectral sequences, the morphism
(\ref{E:mor-spec-seq}) commutes with the differentials in spectral
sequences. We will apply this to the second differential
$d_{\otimes, 2}(\ome\otimes\Delta \eta)$.

First observe that $d_2(\Delta\eta)$, and hence $d_{\otimes,
2}(\ome\otimes\Delta \eta)$, is well defined, in other words
$d_1^{l,0}(\Delta\eta)=0$. Indeed
$d_1^{l,0}(\Delta\eta)=d_1^{l,0}(d_2\eta)$, and, by the definition
of spectral sequence, $d_2\eta$ takes values in the cohomology space
of $d_1^{l,0}$ which in our case is equal to the kernel of
$d_1^{l,0}$, namely $d_1^{l,0}\circ d_2=0$.

Thus by the definition of $d_{\otimes, 2}$ we get
$$d_{\otimes, 2}(\ome\otimes
\Delta\eta)=d_2^{k-2,1}\ome\otimes\Delta\eta\pm \ome\otimes
d_2^{l,0}(\Delta\eta).$$ But $d_2^{l,0}\equiv 0$. Hence
$$d_{\otimes, 2}(\ome\otimes
\Delta\eta)=d_2^{k-2,1}\ome\otimes\Delta\eta=\Delta\ome\otimes
\Delta\eta.$$ The last expression is mapped to
$\Delta\ome\wedge\Delta\eta$ under the morphism
(\ref{E:mor-spec-seq}). \qed

\hfill

We will need few more differential operators on $M$ depending only
on the quaternionic structure. They are defined as dual operators to
$\Delta$ with respect to a very general notion of duality which we
remind now.

For the moment we assume that $M$ is a smooth oriented manifold
without any additional structure. Let $\cs$ and $\ct$ be two finite
dimensional vector bundles over $M$, both either real or complex
simultaneously. Let
$$D\colon C^\infty(M,\cs)\to C^\infty(M,\ct)$$
be a linear differential operator with $C^\infty$-smooth
coefficients. Let us denote by $\ome_M$ the complex line bundle of
differential forms of top degree (either real or complex valued,
depending whether $\ct$ and $\cs$ are real or complex). Let us
consider the operator $D^*$ dual to $D$:
$$D^*\colon (C^\infty_c(M,\ct))^*\to (C^\infty_c(M,\cs))^*$$
where the subscript $c$ stays for compactly supported sections.

Next we have the canonical map
$$C^\infty(M,\ct^*\otimes\ome_M)\to (C^\infty_c(M,\ct))^*$$
given by $f\mapsto [\phi\mapsto \int_M<f,\phi>]$. This map is
injective and has dense image in the weak topology. Thus we will
just identify the image of this map with the source space itself:
\begin{eqnarray*}
C^\infty(M,\ct^*\otimes\ome_M)\subset (C^\infty_c(M,\ct))^*,\\
C^\infty(M,\cs^*\otimes\ome_M)\subset (C^\infty_c(M,\cs))^*.
\end{eqnarray*}
It is easy to see that $D^*$ preserves the class of
$C^\infty$-smooth sections. Actually
$$D^*\colon C^\infty(M,\ct^*\otimes \ome_M)\to
C^\infty(M,\cs^*\otimes\ome_M)$$ is a differential operator of the
same order as $D$ with $C^\infty$-smooth coefficients.

\hfill

Let us apply this construction to our quaternionic manifold $M$ and
the Baston operator
$$\Delta\colon C^\infty(M,\wedge^{k-2}\ce_0^*[-k+1])\to
C^\infty(M,\wedge^k\ce_0^*[-k]),\, 2\leq k\leq 2n.$$

Then we get a second order differential operator which will be used
later
$$\Delta^*\colon C^\infty(M,\wedge^k\ce_0[k]\otimes\ome_M)\to
C^\infty(M,\wedge^{k-2}\ce_0[k-1]\otimes\ome_M).$$

Since the line bundle $\ome_M$ has a canonical real structure and
orientation, the bundles $\wedge^p\ce_0[i]\otimes\ome_M$ are
equipped with a real structure for even $p$, and with a quaternionic
structure for odd $p$. The operator $\Delta^*$ commutes with these
structures, since $\Delta$ does. Notice also that
$$\ome_M\simeq (\det\ce^*_0)^{\otimes 2}[-2n].$$

Thus by the definition of the dual operator we have
\begin{eqnarray}\label{E:baston-dual}
\int_M<f,\Delta\xi>=\int_M<\Delta^*f,\xi>
\end{eqnarray}
for any $\xi\in C^\infty(M,\wedge^{k-2}\ce_0^*[-k+1]),\, f\in
C^\infty(M,\wedge^k\ce_0[k]\otimes\ome_M)$.

\section{Positive currents on quaternionic
manifolds.}\label{S:positive-currents} Let $M^{4n}$, $n>1$, be a
quaternionic manifold. We introduce in this section the notion of
positive (generalized) sections of the bundles
$\wedge^{2k}\ce_0^*[-2k]$, $0\leq k\leq n$. This is analogous to the
notion of positive current from the complex analysis. The case $k=1$
will be necessary for the definition of quaternionic
plurisubharmonic function. In fact the other $k$'s will be important
too, e.g. in the statement (and the proof) of the Aleksandrov and
the Chern-Levine-Nirenberg type theorems.

\hfill

Most of the discussion is actually purely linear algebraic. Let $E$
be a right $\HH$-module of rank $n$, and let $H$ be a left
$\HH$-module of rank 1. As in Section \ref{Ss:quatenionic_mflds} we
will consider $E$ and $H$ as $\CC$-vector spaces via the imbedding
of fields $\CC\inj\HH$ given by $\sqrt{-1}\mapsto I$. Denote
$$E^*:=Hom_\RR(E,\RR),\, H^*:=Hom_\RR(H,\RR);$$
they are left and right $\HH$-modules respectively.

Let $0\leq k\leq n$. The space
$\wedge^{2k}E^*[-2k]:=\wedge_\CC^{2k}E^*\otimes_\CC(\det
H^*)^{\otimes 2k}$ has the real structure defined as follows
(compare with the discussion at the end of Section
\ref{Ss:quatenionic_mflds}). Let $\rho_E$ be the operator on $E$ of
the right multiplication by $j\in \HH$, and let $\rho_H$ be the
operator on $H$ of the right multiplication by $j\in \HH$. Then
$\wedge^{2k}\rho_E^*\otimes (\wedge^{2}\rho_H^*)^{\otimes 2k}$ is
the real structure on $\wedge^{2k}E^*[-2k]$.

The subspace of real elements with respect to this real structure
will be denoted by $\wedge^{2k}E^*[-2k]_\RR$. In this space
$\wedge^{2k}E^*[-2k]_\RR$ we are going to define convex cones
$K^k(E)$ and $C^k(E)$ of weakly positive and strongly positive
elements respectively (the notation might look a bit misleading:
these sets depend of course on the space $H$ too). The cones satisfy
$$C^k(E)\subset K^k(E)\subset \wedge^{2k}E^*[-2k]_\RR,$$
and for $k=0,1,n-1,n$, $C^k(E)=K^k(E)$. The cones are essentially
the same as in \cite{alesker-adv-05}, Section 2.2, though the
construction presented here is simpler.

The definition in the case $k=0$ is obvious: in this case
$\wedge^{2k}E^*[-2k]_\RR=\RR$ and the positive elements are the
usual ones.

Consider the other easy case $k=n$. Clearly $\dim_\CC
\wedge^{2n}E^*[-2n]=1$. Let us fix an $\HH$-bases $e_1,\dots,e_n$ in
$E^*$ and $h$ in $H^*$. Then one can easily see that
$\left(\wedge_{i=1}^n(e_i\wedge Je_i)\right)\otimes (h\wedge
hJ)^{\otimes 2n}$ is a real element of $\wedge^{2n}E^*[-2n]$ and
spans $\wedge^{2n}E^*[-2n]_\RR$. We will call this element positive.
We define the cones $C^n(E)=K^n(E)$ to be the half-line of
non-negative multiples of this element. It is easy to see that this
half-line is independent of choice of bases.

\hfill

Next assume that $1\leq k\leq n-1$. Let us notice that an
$\HH$-linear map $f\colon E\to U$ to another right $\HH$-module
induces a $\CC$-linear map
$$f^*\colon \wedge^{2k}U^*[-2k]\to \wedge^{2k}E^*[-2k]$$
which preserves the real structure. (More precisely $f^*$ is defined
to be $\wedge^{2k}_\CC f^*\otimes Id_{(\det H^*)^{\otimes 2k}}$.)

\begin{definition}\label{D:positive-lin-alg}
(1) An element $\eta\in \wedge^{2k}E^*[-2k]_\RR$, $1\leq k\leq n-1$,
is called {\itshape strongly positive} is it can be presented as a
finite sum of elements of the form $f^*(\xi)$ where $f\colon E\to U$
is a morphism of right $\HH$-modules, $\dim_\HH U=k$, and $\xi\in
\wedge^{2k}U^*[-2k]_\RR$, $\xi\in C^k(U)=K^k(U)$. (Notice that in
particular $\eta=0$ is strongly positive.)

(2) An element $\eta\in \wedge^{2k}E^*[-2k]_\RR$ is called {\itshape
weakly positive} if for any strongly positive
$\xi\in\wedge^{2(n-k)}E^*[-2(n-k)]_\RR$ the wedge product
$\eta\wedge\xi$ is strongly (=weakly) positive element of
$\wedge^{2n}E^*[-2n]_\RR$.
\end{definition}

We have the following properties which are proved in
\cite{alesker-adv-05}, Section 2.2.
\begin{proposition}\label{P:cones-properties}
(1) The cones $C^k(E)$ and $K^k(E)$ are $Aut_\HH(E)\times
Aut_\HH(H)\simeq GL_n(\HH)\times GL_1(\HH)$-invariant, both have
non-empty interior in $\wedge^{2k}E^*[-2k]$, and their closures
contain no $\RR$-linear non-zero subspaces (notice also that
$K^k(E)$ is closed).

(2) $C^k(E)\subset K^k(E)$.

(3) $C^k(E)\wedge C^l(E)\subset C^{k+l}(E)$.

(4) For $k=0,1,n-1,n$
$$C^k(E)=K^k(E).$$
\end{proposition}

For a quaternionic manifold $M$ we denote naturally by
$\wedge^{2k}\ce^*_0[-2k]_\RR$ the vector bundle over $M$ whose fiber
over $p\in M$ is equal to $\left((\wedge^{2k}\ce^*_0|_p)\otimes
(\det\ch_0^*|_p)^{\otimes 2k}\right)_\RR$.

\begin{definition}\label{D:positive-cont}
A continuous section of $\wedge^{2k}\ce^*_0[-2k]_\RR$ is weakly
(resp. strongly) positive if its value at every point is weakly
(resp. strongly) positive.
\end{definition}

\hfill

Let us discuss now positive currents. The discussion is analogous to
the complex case. First let us remind the definition of generalized
section of a vector bundle. Let us discuss this only in the case of
complex vector bundles; the real case is practically the same. Let
$M$ be a smooth manifold which we will assume to be oriented for
simplicity. Let $\cl$ be a (finite dimensional) complex vector
bundle over $M$. Let $C^\infty_c(M,\cl)$ denote the space of smooth
compactly supported sections of $\cl$. This space has a natural
standard linear topology of inductive limit of Fr\'echet spaces. One
denotes by $\ome_M$ the complex line bundle of complex valued
differential forms of top degree. One denotes
$$C^{-\infty}(M,\cl):=(C^\infty_c(M,\cl^*\otimes\ome_M))^*.$$
be the continuous dual. As in Section \ref{S:baston-mult} we have
the natural injective imbedding
$$C^\infty(M,\cl)\inj C^{-\infty}(M,\cl)$$
given by $f\mapsto[\phi\mapsto \int_M <f,\phi>]$. Elements of
$C^{-\infty}(M,\cl)$ are called {\itshape generalized sections} of
$\cl$. The image of this map is dense in the weak topology in the
space of generalized sections.

\hfill

Let us return back to a quaternionic manifold $M$. One has
$$\ome_M\simeq (\det\ce^*_0)^{\otimes 2}[-2n].$$
The wedge product gives a linear map of vector bundles
$$\wedge^{2k}\ce^*_0[-2k]\otimes_\CC
\wedge^{2(n-k)}\ce^*_0[-2(n-k)]\to \wedge^{2n}\ce^*_0[-2n]=(\det
\ce^*_0)[-2n].$$ By duality this map induces a map on vector bundles
which is an isomorphism
$$(\wedge^{2k}\ce^*_0[-2k])^*\tilde\to
\wedge^{2(n-k)}\ce^*_0[-2(n-k)]\otimes_\CC\det\ce_0[2n]=
\wedge^{2(n-k)}\ce^*_0\otimes_\CC(\det\ce_0)[2k].$$ Hence we get an
isomorphism
\begin{eqnarray}\label{E:nov10-1}
(\wedge^{2k}\ce^*_0[-k])^*\otimes_\CC\ome_M\simeq
\wedge^{2(n-k)}\ce^*_0[-2(n-k)]\otimes_\CC\det\ce^*_0.
\end{eqnarray}
All these spaces are equipped with the real structures which are
preserved under the isomorphism (\ref{E:nov10-1}). The real line
bundle of real elements of $\det\ce^*_0$ is canonically oriented.
Hence we can define the convex cones of weakly  (resp. strongly)
positive elements in fibers of
$\wedge^{2(n-k)}\ce^*_0[-2(n-k)]\otimes_\CC\det\ce^*_0$ by taking
tensor products of weakly (resp. strongly) positive elements in
fibers of $\wedge^{2(n-k)}\ce^*_0[-2(n-k)]$ (in the sense of
Definition \ref{D:positive-lin-alg}) with a positive generator of
$\det \ce^*_0$. Via the isomorphism (\ref{E:nov10-1}) this defines
the convex cones of weakly and strongly positive elements in fibers
of $(\wedge^{2k}\ce^*_0[-2k])^*\otimes_\CC\ome_M$. A continuous
section of $(\wedge^{2k}\ce^*_0[-2k])^*\otimes_\CC\ome_M$ is called
weakly (resp. strongly) positive if it is weakly (resp. strongly)
positive at every point.

\begin{definition}\label{D:positive-currents}
A generalized section $\xi\in
C^{-\infty}(M,\wedge^{2k}\ce_0^*[-2k])$ is called {\itshape weakly
positive} (or just {\itshape positive}) if for any strongly positive
smooth compactly supported section $\phi\in
C^\infty_c(M,(\wedge^{2k}\ce_0^*[-2k])^*\otimes\ome_M)$ one has
$\xi(\phi)\geq 0$. We also write in this case: $\xi\geq 0$.
\end{definition}
\begin{remark}\label{R:positive-smooth-gener}
It is easy to see that if $\xi$ is a continuous section, then the
positivity of $\xi$ in the sense of Definition
\ref{D:positive-currents} is equivalent to the weak positivity of
$\xi$ in the sense of Definition \ref{D:positive-cont}.
\end{remark}

\section{Pluripotential theory.}\label{S:psh}
In this section we introduce plurisubharmonic functions on a
quaternionic manifold $M^{4n},\, n>1$. To be more precise these are
not functions but sections of a real line bundle over $M$ defined as
follows.

Let us denote by $(\det \ch_0^*)_\RR$ the real line bundle of real
elements in the bundle $\det\ch_0^*$ which was discussed at the end
of Section \ref{Ss:quatenionic_mflds}. From now on we denote by
$\Delta$ the Baston operator multiplied by an appropriate
normalizing constant (to be chosen in Section \ref{S:hypercomplex}
below) in order to satisfy the compatibility conventions for the
flat manifolds.

\begin{definition}\label{D:psh}
(i) A continuous section $f\in C(M,(\det\ch_0^*)_\RR)$ is called
{\itshape plurisubharmonic} if for the Baston operator $\Delta f\in
C^{-\infty}(M,\wedge^2\ce_0^*[-2]_\RR)$ is positive in the sense of
Definition \ref{D:positive-currents}: $\Delta f\geq 0$.

(ii) A $C^2$-smooth section $f\in C^2(M,(\det\ch^*_0)_\RR)$ is
called {\itshape strictly plurisubharmonic} if at every point $x\in
M$, $\Delta f(x)$ belongs to the interior of the cone of strongly
(=weakly) positive elements of $\wedge^2\ce_0^*[-2]_\RR$.

(iii) A generalized section $f$ of $(\det\ch^*_0)_\RR$ is called
{\itshape pluriharmonic} if $\Delta f=0$.
\end{definition}

\begin{remark}
One can show that any pluriharmonic generalized section $f$ of
$(\det\ch^*_0)_\RR$ is infinitely smooth. This follows from general
elliptic regularity results (see e.g. \cite{aubin}, Ch. 3, \S 6.2,
Theorem 3.54).
\end{remark}

In order to state the Chern-Levine-Nirenberg type estimate let us
fix an auxiliary Riemannian metric on $M$ and metrics on the bundles
$\ce_0,\ch_0$ (or on their relevant tensor products in case
$\ce_0,\ch_0$ are not defined globally).

For any vector bundle $\cl$ with a fixed metric we have norms on the
spaces of (say, continuous) sections of $\cl$:
\begin{eqnarray*}
||\phi||_{L^1(M)}:=\int_M|\phi(x)|dvol(x),\\
||\phi||_{L^\infty(M)}:=\sup_{x\in M}|\phi(x)|.
\end{eqnarray*}
\begin{proposition}[Chern-Levine-Nirenberg type
estimate]\label{P:CLN-estimate} Let $M^{4n}$, $n>1$, be a
quaternionic manifold. Let $1\leq k\leq n$. Let $K, L$ be compact
subsets of $M$ such that $K$ is contained in the interior of $L$.
Then there exists a constant $C$ depending only on auxiliary metrics
and $K,L$ such that for any $C^2$-smooth plurisubharmonic sections
$f_1,\dots,f_k\in C^2(M,(\det\ch^*_0)_\RR)$ one has
$$||\Delta f_1\wedge\dots\wedge \Delta f_k||_{L^1(K)}\leq
C\prod_{i=1}^k||f_i||_{L^\infty(L)}.$$
\end{proposition}
{\bf Proof.} We prove the statement by the induction in $k$. The
case $k=0$ is trivial. Now let us assume the result for $k$
functions and let us prove it for $k+1$. Let us fix a compact subset
$K_1$ containing $K$ and contained in the interior of $L$. Let us
choose $\gamma\in
C^\infty\left(M,(\wedge^{2(k+1)}\ce^*_0[-2(k+1)])^*\otimes
\ome_M\right)$ to be strongly positive (as in Section
\ref{S:positive-currents}), supported in $K_1$, and at very point
$x\in K$ the value $\gamma(x)$ belongs to the interior of the cone
of strongly positive elements. Then there exists a constant $C_1$
such that for any weakly positive continuous section $\xi\in
C(M,\wedge^{2(k+1)}\ce_0^*[-2(k+1)])$ one has
$$||\xi||_{L^1(K)}\leq C_1\int_K<\gamma,\xi>.$$
For any continuous plurisubharmonic sections $f_1,\dots,f_{k+1}$ we
have
\begin{eqnarray*}
||\Delta f_1\wedge\dots \wedge \Delta f_{k+1}||_{L^1(K)}\leq
C_1\int_K<\gamma, \Delta f_1\wedge\dots\wedge \Delta f_{k+1}>\leq\\
C_1\int_{K_1} <\gamma,\Delta f_1\wedge \dots\wedge \Delta
f_{k+1}>\overset{Prop. \, \ref{P:baston-mult}}{=}
C_1\int_{K_1}<\gamma,\Delta(\Delta f_1\wedge\dots\wedge \Delta
f_k\wedge f_{k+1})>=\\
C_1\int_{K_1}<\Delta^*\gamma,\Delta f_1\wedge\dots\wedge \Delta
f_k\wedge f_{k+1}>\leq C_2||f_{k+1}||_{L^\infty(K_1)}\cdot ||\Delta
f_1\wedge\dots\wedge \Delta f_k||_{L^1(K_1)}
\end{eqnarray*}
where $C_2$ depends on $C_1$ and
$||\Delta^*\gamma||_{L^\infty(K_1)}$. Notice that the expression
$\Delta(\Delta f_1\wedge\dots\wedge \Delta f_k\wedge f_{k+1})$ in
the second line is understood in the generalized sense, i.e. as a
$C^{-\infty}$-section of an appropriate vector bundle. Now the rest
follows by the induction assumption. \qed

\hfill

Generalizing the notation from \cite{alesker-verbitsky-06} in the
hypercomplex case, we denote by $P'(M)$ the set of all continuous
plurisubharmonic sections of $(\det \ch_0^*)_\RR$. We denote by
$P''(M)$ the subset of $P'(M)$ consisting of continuous
plurisubharmonic sections $h$ with the following additional
property: for any point $x\in M$ there exist a neighborhood $U$ and
a sequence $\{h_N\}$ of $C^2$-smooth strictly plurisubharmonic
sections over $U$ such that $\{h_N\}$ converges to $h$ uniformly on
compact subsets of $U$ (i.e. in the $C^0$-topology).
\begin{remark}
(1) If $M$ is a locally flat quaternionic manifold, i.e. locally
isomorphic to $\HH^n$, then $P'(M)=P''(M)$. This can be proved
easily by considering convolutions with smooth non-negative
functions.

(2) It is natural to expect that $P'(M)=P''(M)$ for any quaternionic
manifold $M$.

(3) It is easy to see that every section $h\in P''(M)$ locally can
be approximated in the $C^0$-topology by $C^\infty$-smooth strictly
plurisubharmonic sections.
\end{remark}

The following result is an analogue of Proposition 7.8 in
\cite{alesker-verbitsky-06} proved in the hypercomplex case.
\begin{proposition}\label{2-8}
Let $\{h_N\}\subset C(M,(\det\ch_0^*)_\RR)$. Let
$h_N\overset{C^0}{\to}h$. Then

(1) if for any $N$, $h_N\in P'(X)$ then $h\in P'(X)$;

(2) if for any $N$, $h_N\in P''(X)$ then $h\in P''(X)$.
\end{proposition}
{\bf Proof.} Part (2) easily follows from part (1). Thus let us
prove part (1). We have to show that $h$ is plurisubharmonic. Let
$\phi\in C^\infty_0(M,(\wedge^2\ce_0^*[-2])^*\otimes\ome_M)$ be an
arbitrary strongly (=weakly) positive section. We have to show that
$<\phi,\Delta h>\geq 0$. We have
\begin{eqnarray*}
<\phi,\Delta h>=\int_M \Delta^*\phi\cdot h=\lim_{N\to
\infty}\int_M\Delta^*\phi\cdot h_N= \lim_{N\to
\infty}\int_M\phi\cdot \Delta h_N\geq 0.
\end{eqnarray*}
\qed


\hfill

Let us introduce one more notation. We denote by $\tilde
C(M,\wedge^{2k}\ce^*_0[-2k])$ the continuous dual space to the space
of {\itshape continuous} compactly supported sections
$C_0(M,(\wedge^{2k}\ce_0^*[-2k])^*\otimes\ome_M)$. The last space is
equipped with the topology of inductive limit of Banach spaces; that
means that a sequence $\{\xi_i\}$ of continuous compactly supported
sections converges in
$C_0(M,(\wedge^{2k}\ce_0^*[-2k])^*\otimes\ome_M)$ to another such
section if and only if all their supports are contained in some
compact subset and $\xi_i\to\xi$ uniformly. Then $\tilde
C(M,\wedge^{2k}\ce_0^*[-2k])$ is equipped with the weak topology. We
have a natural continuous map
$$\tilde C(M,\wedge^{2k}\ce_0^*[-2k])\to
C^{-\infty}(M,\wedge^{2k}\ce_0^*[-2k])$$ which is injective and has
dense image. We will identify
$$\tilde C(M,\wedge^{2k}\ce_0^*[-k])\subset
C^{-\infty}(M,\wedge^{2k}\ce_0^*[-k]).$$ The following lemma is
essentially well known even in a greater generality, see e.g.
Proposition 5.4 in \cite{alesker-verbitsky-06}.
\begin{lemma}
If $\xi\in C^{-\infty}(M,\wedge^{2k}\ce_0^*[-2k])$ is positive then
it belongs to $\tilde C(M,\wedge^{2k}\ce_0^*[-2k])$.
\end{lemma}

A section $\xi\in \tilde C(M,\wedge^{2k}\ce_0^*[-2k])$ is called
{\itshape positive} if its image in
$C^{-\infty}(M,\wedge^{2k}\ce_0^*[-2k])$ is positive; or
equivalently for any continuous strongly positive $\phi\in
C_0(M,(\wedge^{2k}\ce_0^*[-2k])^*\otimes\ome_M)$ one has
$<\xi,\phi>\geq 0$.


\begin{theorem}\label{T:a-cln}
Let $1\leq k\leq n$. For any $h_1,\dots,h_k\in P''(M)$ one can
define a positive element denoted by $\Delta h_1\wedge
\dots\wedge\Delta h_k\in \tilde C(M,\wedge^{2k}\ce_0^*[-2k]_\RR)$
which is uniquely characterized by the following two properties:

(1) if $h_1,\dots,h_k$ are $C^2$-smooth, then the definition is
straightforward, i.e. pointwise;

(2) if $\{h_i^N\}\subset C^2(M,\wedge^2\ce_0^*[-2]_\RR)$ are
plurisubharmonic, $h_i^N\overset{C^0}{\to}h_i$ as $N\to \infty$,
$i=1,\dots, k$, then
\begin{eqnarray}\label{E:limit-a-cln} \Delta
h_1^N\wedge \dots\wedge \Delta h_k^N\to \Delta h_1\wedge \dots\wedge
\Delta h_k, \, N\to\infty,\end{eqnarray} in the weak topology on
$\tilde C(M,\wedge^{2k}\ce_0^*[-2k]_\RR)$.
\end{theorem}
{\bf Proof.} Let $h_1,\dots,h_k\in P''(M)$. Replacing $M$ by a
smaller open subset if necessary let us choose sequences
$$\{h_i^N\}_{N=1}^\infty\in P'(M)\cap C^2(M,\det\ch_0^*), \,
i=1,\dots,k,$$ such that $h_i^N\overset{C^0}{\to} h_i$ for any
$i=1,\dots,k$. Let us show that $\prod_{i=1}^k\Delta h_i^N$
converges weakly in $\tilde C(M,\wedge^{2k}\ce_0^*[-2k])$ to a
positive element $\mu\in\tilde C(M,\wedge^{2k}\ce_0^*[-2k]_\RR)$.
Since the sequence $\{\prod_{i=1}^k\Delta h_i^N\}_{N=1}^\infty$ is
positive and locally bounded in $L^1$ by Proposition
\ref{P:CLN-estimate}, there exists a subsequence $\{N_l\}$ which
converges weakly to a positive element $\mu\in\tilde
C(M,\wedge^{2k}\ce_0^*[-2k]_\RR)$, and the weak convergence is
understood in the sense of the last space. This fact is a
straightforward generalization of the classical fact that the set of
non-negative measures with bounded integral is compact in the weak
topology. Let us show that $\mu$ does not depend on a choice of
convergent subsequence. This will be shown by induction in $k$. For
$k=0$ this is obvious. Let us assume the statement for $k-1$ and
prove for $k$. Let $\nu$ be a weak limit of another subsequence. It
suffices to check that for any $\phi\in
C^\infty_0(M,(\wedge^{2k}\ce_0^*[-2k])^*\otimes\ome_M)$ one has
\begin{eqnarray}\label{E:nov12-1}
<\mu,\phi>=<\nu,\phi>.
\end{eqnarray}
We have
\begin{eqnarray}
<\mu,\phi>=\lim_{l\to\infty}\int_M<\prod_{i=1}^k\Delta
h_i^{N_l},\phi>\overset{\mbox{Prop. }
\ref{P:baston-mult}}{=}\\\label{Delta_Aug11} \lim_{l\to
\infty}\int_M<\Delta\left((\prod_{i=1}^{k-1} \Delta h_i^{N_l})\wedge
h_k^{N_l}\right),\phi>=\\\label{E:nov12-2} \lim_{l\to
\infty}\int_M<\prod_{i=1}^{k-1}\Delta h_i^{N_l}, h_k^{N_l}\cdot
\Delta^* \phi>.
\end{eqnarray}
Notice that in (\ref{Delta_Aug11}) the expression
$\Delta\left((\prod_{i=1}^{k-1} \Delta h_i^{N_l})\wedge
h_k^{N_l}\right)$ is understood in the generalized sense, i.e. as a
$C^{-\infty}$-section of an appropriate vector bundle.

By the induction assumption the sequence
$g_N:=\prod_{i=1}^{k-1}\Delta h_i^{N}$ weakly converges to an
element $g$. The sequence $\{f_N:=h_k^N\wedge \Delta^*\phi\}$ has
uniformly bounded support and converges uniformly (i.e. in the
$C^0$-topology) to $f:=h_k\wedge \Delta^*\phi$. Thus existence of
the limit (\ref{E:nov12-2}) follows from the following well known
lemma (see e.g. Lemma 7.12 in \cite{alesker-verbitsky-06}).
\begin{lemma}\label{L:nov12-3}
Let $X$ be a compact topological space. Let $E\to X$ be a finite
dimensional vector bundle. Let $\{f_N\}\subset C(X,E)$ be a sequence
of continuous sections which converges to $f\in C(X,E)$ uniformly on
$X$. Let $\{g_N\}\subset C(X,E)^*$ be a sequence in the dual space
which weakly converges to $g\in C(X,E)^*$. Then $<g_N,f_N>\to <g,f>$
as $N\to\infty$.
\end{lemma}
Let us only notice that this lemma is an easy consequence of the
Banach-Steinhauss theorem. This lemma implies that the limit
$\lim_{N\to\infty}\int_M<\prod_{i=1}^{k-1}\Delta h_i^N, h_k^N\wedge
\Delta^*\phi>$ does exist, and the same argument shows that it is
equal also to $<\nu,\phi>$. Hence the equality (\ref{E:nov12-1}) is
proved. Hence there exists a weak limit $\mu$ of the sequence
$\prod_{i=1}^k\Delta h_i^N$.

\hfill

Let us show that if $h_1,\dots,h_k$ are $C^2$-smooth then the limit
$\mu$ is equal to $\prod_{i=1}^k\Delta h_i$ defined just pointwise.
The proof is again by induction on $k$. The case $k=0$ is trivial.
To make the induction step, let us fix $\phi\in
C^\infty_0(M,(\wedge^{2k}\ce_0^*[-2k])^*\otimes\ome_M)$. We have
\begin{eqnarray}\label{E:nov12-10}
<\mu,\phi>=\lim_{N\to \infty}\int<\prod_{i=1}^k\Delta h_i^N,\phi>\overset{\mbox{ Prop. } \ref{P:baston-mult}}{=}\\
\lim_{N\to \infty}\int<\Delta\left((\prod_{i=1}^{k-1}\Delta
h_i^N)\wedge h_k^N\right),\phi>=\\
\lim_{N\to \infty}\int<\prod_{i=1}^{k-1}\Delta h_i^N, h_k^N\cdot
\Delta^*\phi>.\label{E:nov12-11}
\end{eqnarray}
By Lemma \ref{L:nov12-3} and the induction assumption the last limit
exists and is equal to
\begin{eqnarray}\label{E:nov12-12}
\int<\prod_{i=1}^{k-1}\Delta h_i, h_k\cdot\Delta^*\phi>.
\end{eqnarray}
Let us show that (\ref{E:nov12-12}) is equal to
$\int<\prod_{i=1}^k\Delta h_i,\phi>$. But this is proved exactly in
the same way as the equality (\ref{E:nov12-10})=(\ref{E:nov12-11}).
To complete the proof of the theorem, it remains to prove the
following result which we would like to formulate as a separate
statement.

\begin{theorem}\label{Tt:convergence}
Let $\{h_i^N\}_{N=1}^\infty\subset P''(M)$, $i=1,\dots,k$, be a
sequence such that for any $i=1,\dots,k$
$$h_i^N \overset{C^0}{\to}h_i.$$
Then $h_i\in P''(M)$ and
$$\prod_{i=1}^k\Delta h_i^N\to \prod_{i=1}^k\Delta h_i$$
in the space $\tilde C(M,\wedge^{2k}\ce_0^*[-2k])$, and the products
are defined in the sense of Theorem \ref{T:a-cln}.
\end{theorem}
This theorem is proved by inspection of the proof of Theorem
\ref{T:a-cln}. \qed

\begin{remark}
Theorems \ref{T:a-cln} and \ref{Tt:convergence} have classical real
and complex analogs for convex and complex plurisubharmonic
functions respectively. The real analogue is due to A.D. Aleksandrov
\cite{alexandrov-58}, and the complex one is due to
Chern-Levine-Nirenberg \cite{chern-levine-nirenberg}. Notice also
that Theorems \ref{T:a-cln} and \ref{Tt:convergence} were proved in
the special case of flat quaternionic space in
\cite{alesker-bsm-03}, and their analogue for hypercomplex manifolds
in \cite{alesker-verbitsky-06}.
\end{remark}

\begin{definition}
Let us define the Monge-Amp\`ere operator on sections of $(\det
\ch_0^*)_\RR$ by
$$h\mapsto (\Delta h)^n.$$
\end{definition}
The Monge-Amp\`ere operator is naturally defined on $C^2$-smooth
sections of $(\det \ch_0^*)_\RR$. But by Theorem \ref{T:a-cln} it
can be defined for sections from $P''(M)$. Using this operator one
can easily introduce a Monge-Amp\`ere equation on any quaternionic
manifold, but we do not pursue this point here.

\section{The case of flat
manifolds.}\label{S:hypercomplex} A basic theory of plurisubharmonic
functions on the flat space $\HH^n$ was developed by the author in
\cite{alesker-bsm-03} (see also \cite{alesker-jga-03},
\cite{alesker-adv-05}), and on more general class of hypercomplex
manifolds by M. Verbitsky and the author
\cite{alesker-verbitsky-06}. The theory of \cite{alesker-bsm-03} is
a special case of \cite{alesker-verbitsky-06}. In this section we
show that the theory of the flat case in \cite{alesker-bsm-03} is a
special case of the theory of this paper too. Since hypercomplex
manifolds form a subclass of quaternionic manifolds, it is natural
to compare the classes of plurisubharmonic functions on hypercomplex
manifolds from \cite{alesker-verbitsky-06} and of the present paper.
If the holonomy of the Obata connection is not contained in
$SL_n(\HH)$ then the theories are formally different since
plurisubharmonic sections belong to different line bundles. For the
moment we do not know whether the two classes do coincide when the
holonomy of the Obata connection is contained in $SL_n(\HH)$. It
would be interesting to give a sufficient condition on a
hypercomplex manifold under which the two theories can be
identified.

We need some preparations.
Let us make some identifications in the special case of $M=\HH^n$.
Recall that by (\ref{qm5})
\begin{eqnarray}\label{E:isom-Jule2011}
TM\simeq \ce_0\otimes_\HH \ch_0.
\end{eqnarray}
Let $X$ denote a complexification of $M=\HH^n$, as previously. As
previously, we denote by $\ce,\ch$ the holomorphic vector bundles
over $X$ extending the real analytic vector bundles $\ce_0,\ch_0$
respectively. $\ce_0,\ch_0$ are topologically trivial quaternionic
vector bundles, while $\ce,\ch$ are holomorphically trivial
(complex) vector bundles. Moreover $\ce_0,\ch_0$ are equivariant
under the group of affine transformations $\ca:=\HH^n\rtimes
(SL_n(\HH)\times SL_1(\HH))$, i.e. the semi-direct product of the
group of translations $\HH^n$ and linear transformations
$SL_n(\HH)\times SL_1(\HH)$ (here $SL_1(\HH)=\{q\in \HH|\,
|q|=1\}$); notice that the factor $SL_1(\HH)$ acts trivially on
$\ce_0$ while the factor $SL_n(\HH)$ acts trivially on $\ch_0$.
Moreover $\ca$ acts by automorphisms of $\HH^n$ as a quaternionic
manifold, and the induced action on the tangent bundle is compatible
with the isomorphism (\ref{E:isom-Jule2011}). Hence all the Baston
operators become equivariant under this group $\ca$. This property
will be crucial in the proof of the main result of this section.
Thus in particular the operator
$$\Delta\colon C^\infty(M,\det\ch_0^*)\to
C^\infty(M,\wedge^2\ce_0^*[-2])$$ is $\ca$-equivariant.

It is easy to see that the complex line bundle $\det\ch_0^*$ is
isomorphic to the trivial bundle in the category of
$\ca$-equivariant vector bundles. Moreover the isomorphism can be
chosen in such a way that the real structure on fibers of
$\det\ch_0^*$ (which is induced by the quaternionic structure on
fibers on $\ch_0$) is mapped to the real structure on the trivial
line bundle, and positive half line is mapped to the positive half
line.

Next the complex vector bundle $\wedge^2\ce_0^*[-2]$ is
$\ca$-equivariantly isomorphic to the vector bundle over $\HH^n$
whose fiber over a point $p\in \HH^n$ is equal to the space of
$\CC$-valued quaternionic Hermitian forms on $T_p\HH^n$
(equivalently, to $\CC$-valued quadratic forms on the real space
$T_p\HH^n$ which are invariant under group of norm one quaternions);
we refer to \cite{alesker-verbitsky-06} for this linear algebra.
Moreover this isomorphism can be chosen to preserve the real
structures on both spaces and the cones of positive elements (the
positive cone in the latter space was defined in
\cite{alesker-bsm-03}, and an equivalent description in the more
general case of hypercomplex manifolds was given in
\cite{alesker-verbitsky-06}; the definitions of positivity in the
present paper just a direct generalization of the definitions from
\cite{alesker-verbitsky-06}). Notice that the quaternionic Hessian
defined in \cite{alesker-bsm-03} in the flat case took values
exactly in this vector bundle for quaternionic Hermitian forms. With
these identifications, the Baston operators $\Delta$ and the
quaternionic Hessian introduced in \cite{alesker-bsm-03} act between
the same vector bundles. We will denote the quaternionic Hessian
from \cite{alesker-bsm-03} by $\Delta'$. By \cite{alesker-bsm-03},
the operator $\Delta'$ is also $\ca$-equivariant.

\begin{proposition}
Let $M=\HH^n$. With the above identifications

(i) the Baston operator $\Delta$ (appropriately normalized)
coincides with $\Delta'$ from \cite{alesker-bsm-03};

(ii) the Monge-Amp\`ere operator $h\mapsto (\Delta h)^n$ coincides
with the Monge-Amp\`ere operator from \cite{alesker-bsm-03};

(iii) the class of plurisubharmonic functions in the sense of
Definition \ref{D:psh} of this paper coincides with the class of
plurisubharmonic functions introduced in \cite{alesker-bsm-03}.
\end{proposition}
{\bf Proof.} Parts (ii),(iii) follow from part (i) and the
definitions of the real structures, positive cones, and the wedge
product (which is equivalent by \cite{alesker-verbitsky-06} to the
Moore determinant of quaternionic matrices used in
\cite{alesker-bsm-03}). Part (i) is the main one, and we are going
to prove it. We have to show that $\Delta'=\Delta$ when $\Delta$ is
appropriately normalized.

\hfill

First let us prove the vanishing of the symbol of $\Delta'-\Delta$
considered as a differential operator of second order. By the
translation invariance it suffices to show that the symbol of
$\Delta'-\Delta$ vanishes at 0. Since $\Delta'-\Delta$ is
$\ca$-equivariant, its symbol belongs to

\begin{eqnarray}\label{E:001}
Hom_{\RR,\cb}\left(Sym_\RR^2(\HH^{n\ast})\otimes
\det(\ch_0^*)_\RR|_0,(\wedge^2\ce^*_0[-2])_\RR|_0\right)
\end{eqnarray}
where we have denoted for brevity $\cb:=SL_n(\HH)\times SL_1(\HH)$,
and $Hom_{\RR,\cb}$ denotes the space of $\RR$-linear maps commuting
with the group $\cb$. Now we are going to show that the space
(\ref{E:001}) is at most one dimensional. It suffices to show the
one dimensionality of the $Hom$ between complexified
representations, namely that
\begin{eqnarray}\label{E:002}
Hom_{\CC, \ccb}\left(Sym^2_\CC((\CC^{2n}\otimes \CC^2)^*)\otimes
\det(\CC^{2*}),\wedge_\CC^2(\CC^{2n})^*\otimes
(\det\CC^{2*})^{\otimes 2}\right)
\end{eqnarray}
is one dimensional. Here $\ccb=SL_{2n}(\CC)\times SL_2(\CC)$ is the
complexification of $\cb$. Next in (\ref{E:002}) the action of
$\ccb$ on $\CC^{2n}$ and $\CC^2$ is as follows:  $SL_{2n}(\CC)$ and
$SL_2(\CC)$ act on $\CC^{2n}$ and respectively $\CC^2$ in the
standard way, $SL_{2n}(\CC)$ acts trivially on $\CC^2$, and
$SL_2(\CC)$ acts trivially on $\CC^{2n}$.

The representation of $SL_{2n}(\CC)$ on $\wedge^2_\CC(\CC^{2n})^*$
is irreducible (see e.g. \cite{goodman-wallach}, Corollary 5.5.3).
The representation of $SL_{2n}(\CC)\times SL_2(\CC)$, and hence of
$\ccb$, in $Sym^2(\CC^{2n}\otimes \CC^2)$ is multiplicity free (see
e.g. \cite{goodman-wallach}, Corollary 5.6.6); hence the
representation of $\ccb$ in $$Sym^2_\CC((\CC^{2n}\otimes
\CC^2)^*)\otimes \det(\CC^{2*})$$ is also multiplicity free. Then
the Schur's lemma implies that the $Hom$-space (\ref{E:002}) is at
most one dimensional. This implies that the symbols of $\Delta'$ and
$\Delta$ must be proportional. Hence $\Delta$ can be normalized in
such a way that the symbols just coincide.

\hfill

Thus the differential operator $\Delta'-\Delta$ has order at most
one. Let us consider the symbol of this first order differential
operator. It is an element of
\begin{eqnarray}\label{E:003}
Hom_{\RR,\cb}\left(\HH^{n*}\otimes
\det(\ch_0^*)_\RR|_0,(\wedge^2\ce_0^*[-2]_\RR)|_0\right).
\end{eqnarray}
But the two representations under the $Hom$ are irreducible and
non-isomorphic. Hence (\ref{E:003}) vanishes, and the symbol of
$\Delta'-\Delta$ vanishes. Hence $\Delta'-\Delta$ has order zero.
But then $\Delta'-\Delta$ defines an element of
$$Hom_{\RR,\cb}\left((\det\ch_0^*)_\RR|_0,\wedge^2\ce^*_0[-2]_\RR|_0\right).$$
Obviously the last space vanishes. This implies that
$\Delta'-\Delta=0$. \qed

\end{document}